\documentclass[11pt]{article}

\renewcommand{\baselinestretch}{1}   %line space adjusted here
\usepackage{psfrag,cite,amsmath,graphicx,epsfig,latexsym,subfigure,color}
\usepackage[scriptsize,bf]{caption}
\usepackage{color,xcolor}
\usepackage{amssymb}
\usepackage{setspace}
\usepackage{array}
\usepackage{booktabs}
\def\bn{\hfill \\ \smallskip\noindent}

\def\dom{\mbox{dom\,}}

\newcommand{\beq}{\begin{equation}}
\newcommand{\eeq}{\end{equation}}
\newcommand{\st}{{\rm s.t.}}

\newcommand{\cG}{{\mbox{$\mathcal{G}$}}}

\newcommand{\cV}{{\mbox{$\mathcal{V}$}}}

\begin{document}
\def\pn {\par\smallskip\noindent}
\def \bn {\hfill \\ \smallskip\noindent}
\newcommand{\fs}{f_1,\ldots,f_s}
\newcommand{\f}{\vec{f}}
\newcommand{\hf}{\hat{f}}
\newcommand{\hx}{\hat{x}}
\newcommand{\hy}{\hat{y}}
\newcommand{\hz}{\hat{z}}
\newcommand{\hw}{\hat{w}}
\newcommand{\tw}{\tilde{w}}
\newcommand{\hlambda}{\hat{\lambda}}
\newcommand{\hbeta}{\hat{\beta}}
\newcommand{\tG}{\widetilde{G}}
\newcommand{\tg}{\widetilde{g}}
\newcommand{\barhx}{\bar{\hat{x}}}
\newcommand{\vecx}{x_1,\ldots,x_m}
\newcommand{\xoy}{x\rightarrow y}
\newcommand{\barx}{{\bar x}}
\newcommand{\bary}{{\bar y}}
\newcommand{\hrho}{\widehat{\rho}}
\newtheorem{theorem}{Theorem}[section]
\newtheorem{lemma}{Lemma}[section]
\newtheorem{corollary}{Corollary}[section]
\newtheorem{proposition}{Proposition}[section]
\newtheorem{definition}{Definition}[section]
\newtheorem{claim}{Claim}[section]
\newtheorem{remark}{Remark}[section]

\def\br{\break}
\def\smskip{\par\vskip 5 pt}
\def\proof{\bn {\bf Proof.} }
\def\QED{\hfill{\bf Q.E.D.}\smskip}
\def\qed{\quad{\bf q.e.d.}\smskip}
\newcommand{\cE}{\mathcal{E}}
\newcommand{\cM}{\mathcal{M}}
\newcommand{\cN}{\mathcal{N}}
\newcommand{\cJ}{\mathcal{J}}
\newcommand{\cT}{\mathcal{T}}
\newcommand{\bx}{\mathbf{x}}
\newcommand{\bp}{\mathbf{p}}
\newcommand{\bX}{\mathbf{X}}
\newcommand{\bY}{\mathbf{Y}}
\newcommand{\bP}{\mathbf{P}}
\newcommand{\bA}{\mathbf{A}}
\newcommand{\bB}{\mathbf{B}}
\newcommand{\bfM}{\mathbf{M}}
\newcommand{\bL}{\mathbf{L}}
\newcommand{\bz}{\mathbf{z}}
\newcommand{\cF}{\mathcal{F}}
\newcommand{\cR}{\mathcal{R}}
\newcommand{\bzero}{\mathbf{0}}

\newcommand{\blue}{\color{blue}}
\newcommand{\red}{\color{red}}

	\title{Gradient Primal-Dual Algorithm Converges to Second-Order Stationary Solutions for  Nonconvex Distributed Optimization}
	\author{{Mingyi Hong \thanks{M. Hong is with the Department of Electrical and Computer Engineering, University of Minnesota, Minneapolis, MN 55414, USA. Email: \texttt{mhong@umn.edu}}, Jason D. Lee \thanks{J. D. Lee is with the Department of Data Sciences and Operations, the University of Southern California, Los Angeles, CA 90089. Email: \texttt{jasonlee@marshall.usc.edu}},  Meisam Razaviyayn\thanks{M. Razaviyayn is with the department of Industrial and Systems Engineering, the University of Southern California. Email: \texttt{razaviya@usc.edu}}}}
	\maketitle

	\begin{abstract}
		In this work, we study two first-order primal-dual based algorithms, the Gradient Primal-Dual Algorithm (GPDA) and the Gradient Alternating Direction Method of Multipliers (GADMM),  for  solving a class of linearly constrained non-convex optimization problems. We show that with random initialization of  the primal and dual variables,  both algorithms are able to compute second-order stationary solutions (ss2) with probability one.
		This is the first result showing that primal-dual algorithm is capable of finding ss2 when only using first-order information;  it also extends the existing results for first-order, but {\it primal-only} algorithms.
		
		An important implication of our result is that it also gives rise to the first global convergence result to the ss2, for two classes of  unconstrained {\it distributed} non-convex learning problems over multi-agent networks. %, a class of primal-dual gradient based algorithm also converges to second-order stationary solutions.
	\end{abstract}
	\vspace{-0.15in}

	\section{Introduction}\label{sec:introduction}
	%\subsection{The problem}
	In this work, we consider the following linearly constrained optimization problem:
	\begin{align}\label{eq:original}
	\min_{x\in\mathbb{R}^N}\; f(x) \quad \st \quad A x=b
	\end{align}
	where $f(x):\mathbb{R}^{N}\to \mathbb{R}$ is a smooth function (possibly non-convex); $A\in\mathbb{R}^{M\times N}$ is not full column rank; $b\in\mathbb{R}^{M}$ is a known vector. %$X$ is a {\it convex} compact set; $h(x):\mathbb{R}^N \to \mathbb{R}$ is a nonsmooth but convex function.
	
	An important application of problem \eqref{eq:original} is in the non-convex distributed optimization  and learning -- a problem that has gained considerable attention recently, and has found applications in  training neural networks \cite{Lian17decentralized}, distributed information processing and machine learning \cite{Forero11,hong14nonconvex_admm}, and distributed signal processing \cite{Lorenzo16}. In distributed optimization and learning, the common setup is that a network consists of $N$ distributed agents collectively optimize the following problem %the following classical global consensus problem
	\begin{align}\label{eq:global:consensus}
	\min_{v\in\mathbb{R}} \quad \sum_{i=1}^{N} f_i(v) + g(v),
	\end{align}
	where $f_i(v): \mathbb{R}\to\mathbb{R}$ is a function local to agent $i$ (note, for notational simplicity we assume that $v$ is a scalar); $g(v)$ represents some smooth regularization function known to all agents. Below we present two problem formulations based on  different topologies and application scenarios.  %{[[[\color{red} I changed $y$ to $v$ to avoid confusion with the next section when we use $y$]]]}
	
	\noindent{\bf Scenario 1: The Global Consensus.} Suppose that all the agents are connected to a single central node. The distributed agents can communicate with the controller, but they are not able to directly communicate among themselves. In this case problem \eqref{eq:global:consensus} can be equivalently formulated into the following global consensus problem  \cite{BoydADMMsurvey2011,hong14nonconvex_admm}
	\begin{align}\label{eq:global:consensus:1}
	\min_{\{x_i\}_{i=0}^{N}}\; \sum_{i=1}^{N}f_i(x_i) + g(x_0), \quad \st \quad x_i =x_0,\; \forall~i.
	\end{align}
	The setting of the above global consensus problem is popular in applications such as parallel computing, in which the existence of central controller can orchestrate the activity of all agents; see \cite{icml2015_zhangb15, li2013distributed}.  To cast the problem into the form of \eqref{eq:original}, define
	\begin{align}\label{eq:global:consensus:equiv:0}
	f(x) & = \sum_{i=1}^{N}f_i(x_i) + g(x_0), \nonumber\\
	&A_1 = I_{N}, \;  A_2 = 1_N, \;  A=[A_1, -A_2], \; b=0,
	\end{align}
	where $I_N\in\mathbb{R}^{N\times N}$ is the identity matrix; $1_N\in\mathbb{R}^{N}$ is the all one vector.

	\noindent{\bf Scenario 2: Distributed Optimization Over Networks.} Suppose that there is no central controller, and the $N$ agents are connected by a network defined by an {\it undirected} graph $\cG=\{\mathcal{V}, \mathcal{E}\}$, with $|\mathcal{V}|=N$ vertices and $|\mathcal{E}|=E$ edges. Each agent can only communicate with its immediate neighbors, and it can access one component function $f_i$. This problem has wide applications ranging from distributed communication networking \cite{liao15semi-async},  distributed and parallel machine learning \cite{Forero11,mateos_dlasso, shalev14proximaldual}, to distributed signal processing \cite{Schizas08}.
	
	Define the node-edge incidence matrix $A\in\mathbb{R}^{E\times N}$ as following: if $e\in\mathcal{E}$ and it connects vertex $i$ and $j$ with $i>j$, then $A_{ev}=1$ if $v=i$, $A_{ev}=-1$ if $v=j$ and $A_{ev}=0$ otherwise.  %Using this definition, the {\it signed graph Laplacian matrix} $L_{-}\in\mathbb{R}^{N\times N}$ is given by
	%$$L_{-} := A^T A.$$
	Introduce $N$ local variables $x=[x_1,\cdots, x_N]^T$, and suppose the graph $\{\cV, \cE\}$ is connected. Then as long as as the graph  is connected, the following formulation is equivalent to the global consensus problem, which is precisely problem \eqref{eq:original}
	\begin{align}\label{eq:global:consensus:equiv}
	\hspace{-0.5cm}\min_{x\in\mathbb{R}^N} \; f(x):=\sum_{i=1}^{N} \left(f_i(x_i)+\frac{1}{N} g(x_i)\right),\; \st\; Ax = 0.
	\end{align}
	
	%{\red [Shall we discuss application in phase retrieval?]}
	
	\subsection{The objective of this work}
	The research question we attempt to address in this work is:
	
	%{\color{red}[[[ Maybe it would be better to say the objective is ``Can we design primal-dual algorithms capable of computing...". This way we exclude projected gradient descent.]]]}
%	\vspace{-0.2cm}
\begin{center}
	\noindent\fcolorbox{black}[rgb]{0.9,0.9,0.9}{\begin{minipage}{5in}
			\begin{center}
				{\bf(Q)} ~Can we design primal-dual algorithms capable of computing second-order stationary solutions for \eqref{eq:original}?
			\end{center}
		\end{minipage}}
\end{center}

	%	\vspace{-0.2cm}

		Let us first analyze the first-order stationary (ss1) and second-order stationary (ss2) solutions for problem \eqref{eq:original}. For a general smooth nonlinear problem in the following form
		\begin{align}\label{eq:nonlinear}
		\min_{x\in\mathbb{R}^N}\; g(x) \quad \st \quad h_i(x)=0, \quad i=1,\cdots, m,
		\end{align}
		the first-order necessary condition is given as
		\begin{align}\label{eq:first_order}
		\hspace{-0.5cm}\nabla g(x^*) + \sum_{i=1}^{m}\langle \lambda_i^*, \nabla h_i(x^*)\rangle = 0, \quad h_i(x^*) = 0, \; \forall~i.
		\end{align}
		The second-order necessary condition is given below [see Proposition 3.1.1 in \cite{bertsekas99}]. Suppose $x^*$ is regular, then
		\begin{align}\label{eq:necessary}
		\begin{split}
		&\langle y, (\nabla^2 g(x^*) +\sum_{i=1}^{m}\lambda^*_i \nabla^2 h_i(x^*))y\rangle \ge 0,\\
		&\forall~y\in \{y\ne 0\mid \langle \nabla h_i(x^*), y\rangle =0, \; \forall~i=1,\cdots, m \}.
		\end{split}
		\end{align}
		
		Applying the above result to our problem, we obtain the following first- and second-order necessary condition for problem \eqref{eq:original} \footnote{Note that for linear constraints no further regularity is needed for the existence of multipliers}
		\begin{subequations}\label{eq:first:and:second}
			\begin{align}
			&\nabla f(x^*) + A^T \lambda^* = 0, \quad A x^* = b. \label{eq:stationarity:ss1}\\
			&\langle y, \nabla^2 f(x^*) y\rangle \ge 0,\quad \forall\;\; y \in \{y\mid A y =0\}.\label{eq:stationarity}
			\end{align}
		\end{subequations}
		In other words, the second-order necessary condition is equivalent to the condition that $\nabla^2 f(x^*)$ is positive semi-definite in the null space of $A$.
		Similarly, the sufficient condition for {\it strict} local minimizer is given by %\cite[Proposition 3.2.1.]{bertsekas99}
		\begin{align}\label{eq:sufficient:1}
		%	\begin{split}
		&\nabla f(x^*) +  A^T \lambda^* = 0, \quad A x^* = b. \\
		&\langle y, \nabla^2 f(x^*) y\rangle > 0,\quad \forall\;\; y\ne 0, \; \mbox{and}\;  y\in \{y\mid A y =0\}.\nonumber
		%	\end{split}
		\end{align}
		%In another word, the second-order necessary condition is equivalent to the condition that $\nabla^2 f(x^*)$ is positive semi-definite in the null space of $A$.
		
		To proceed, we need the following claim [see Lemma 3.2.1 in \cite {bertsekas99}]
		\begin{claim}\label{claim:null:space}
			Let $P$ and $Q$ be two symmetric matrices. Assume that $Q$ is positive semidefinite and $P$ is positive definite on the null space of $Q$, that is, $x^T P x >0$ for all $x\ne 0$ with $x^T Q x =0$. Then there exists a scalar $\bar{c}$ such that
			\begin{align}\label{eq:pd:1}
			P+cQ \succ 0, \quad \forall ~ c\ge \bar{c}.
			\end{align}
			Conversely, if there exists a scalar $\bar{c}$ such that \eqref{eq:pd:1} is true,
			then we have $x^T P x >0$ for all $x\ne 0$ with $x^T Q x =0$.
		\end{claim}

		By Claim \ref{claim:null:space}, the sufficient condition \eqref{eq:sufficient:1}  can be equivalently written as:  %the following condition is {\it sufficient} for $x^*$ being a strong local min:
		\begin{align}\label{eq:sufficient}
		&\nabla f(x^*) + A^T \lambda^* = 0, \quad A x^* = b. \\
		&\nabla^2 f(x^*) + \gamma A^T A \succ 0, \; \mbox{for~some~}\gamma>0.
		\end{align}
		%Therefore, checking the sufficient condition \eqref{eq:sufficient:1} is equivalent to checking the positive definiteness of certain matrices, therefore it is polynomial time; also see
		
		%By Claim \ref{claim:null:space:reverse}, the conditions below implies the necessary condition \eqref{eq:stationarity}
		%\begin{subequations}
		%\begin{align}
		%&\nabla f(x^*) + \langle \lambda^*, A\rangle = 0, \quad A x^* = b. \\
		%&\nabla^2 f(x^*) + \gamma A^T A \succeq 0, \; \mbox{for~some~}\gamma>0.  \label{eq:stationarity:2}
		%\end{align}
		%\end{subequations}
		%If further $\nabla^2 f(x^*)$ is full rank, then \eqref{eq:stationarity:2} and \eqref{eq:stationarity} are equivalent.

		It is worth mentioning that checking both of the above sufficient and necessary conditions can be done in polynomial time, but when there are inequality constraints, checking second-order conditions can be NP-hard; see \cite{Murty1987}.
		In the following we will refer to the condition \eqref{eq:stationarity:ss1} as ss1 solution and condition \eqref{eq:stationarity} as the ss2 solution.
		According to the above definition, we define a {\it strict saddle} point to be the solution $x^*$ such that
		\begin{align}\label{eq:strict}
		\begin{split}
		&\nabla f(x^*) +  A^T \lambda^* = 0, \quad A x^* = b, \\
		&\exists~ y\in\{y\mid Ay=0, y\ne 0\}, \; \mbox{and}\; \sigma>0\quad \mbox{such that}\; \langle y, \nabla^2 f(x^*) y\rangle  \le -\sigma \|y\|^2.
		\end{split}
		\end{align}
		It is easy to verify using Claim \ref{claim:null:space} that the above condition implies that for the same $\sigma>0$, the following is true
		\begin{align}\label{eq:strict:2}
		\begin{split}
		&\nabla f(x^*) +  A^T \lambda^* = 0, \quad A x^* = b, \\
		&\sigma_{\min}\left(\gamma A^T A + \nabla^2 f(x^*) \right)  \le -\sigma, \quad \forall~\gamma>0
		\end{split}
		\end{align}
		%In words, for any $\gamma>0$, the matrix $\nabla^2 f(x^*) + \gamma A^T A$ has a least one negative eigenvalue bounded away from zero.
		where $\sigma_{\min}$ denotes the smallest eigenvalue of a matrix.
		Clearly, if a ss1 solution $x^*$ does not satisfy \eqref{eq:strict}, i.e.,
		\begin{align}
		\forall~y, \; \mbox{s.t.} \; A y =0, \quad \langle y, \nabla f^2(x^*) y\rangle \ge 0,
		\end{align}
		then \eqref{eq:stationarity} is true. In this work, we will develop primal-dual algorithms that avoid converging to the strict saddles \eqref{eq:strict}.
		
		\subsection{Existing literature}
		%\subsubsection{Methods for finding second-order stationary solutions}
		
		Many recent works have been focused on designing algorithms
		with convergence guarantees to local minimum points/ss2 for non-convex
		unconstrained problems. These include second-order methods such as  trust region method
		\cite{conn2000trust}, cubic regularized Newton's method
		\cite{nesterov2006cubic}, and a hybrid of first-order and second-order methods \cite{sama17}. When only gradient information is available, it
		has been shown that with random initialization, gradient descent (GD) converges to ss2 for unconstrained smooth problems with probability one
		\cite{jlee16jordan}. Recently,
		a perturbed version of GD which occasionally adds noise to the iterates has been proposed \cite{jin2017jordan}, and such a method converges
		to the ss2 with faster
		convergence rate than the ordinary gradient descent algorithm with random initialization. When manifold constraints are present, it is shown in \cite{lee17first_order} that manifold gradient descent converges to ss2, provided that each time the iterates are always feasible (ensured by performing a potentially expensive second-order retraction operation).
		However, there has been no work analyzing whether classical primal-dual gradient type methods based on Lagrangian relaxation are also capable of computing ss2.
		
		%\subsubsection{Methods for distributed non-convex optimization}
		
		The consensus problem \eqref{eq:global:consensus} and \eqref{eq:global:consensus:equiv} have been studied extensively in the literature when the objective functions are all convex; see for example \cite{Nedic09subgradient, nedic2015distributed, shi2014extra, aybat2016primal}. Primal methods such as distributed subgradient method \cite{Nedic09subgradient}, the EXTRA method \cite{shi2014extra}, as well as primal-dual based methods such as Alternating Direction Method of Multipliers (ADMM) \cite{BoydADMMsurvey2011,Schizas09,chang14distributed} have been studied.  On the contrary, only recently there have been some work addressing the more challenging problems without assuming convexity of  $f_i$'s; see recent developments in \cite{bianchi2013convergence, hong14nonconvex_admm,Lorenzo16,hong17icml}. In particular, reference \cite{hong14nonconvex_admm} develops non-convex ADMM based methods (with global sublinear convergence rate) for solving the global consensus problem \eqref{eq:global:consensus:1}. Reference \cite{hong17icml} proposes a primal-dual based method for unconstrained non-convex distributed optimization over a connected network (without a central controller), and derives the first global convergence rate for distributed non-convex optimization. In \cite{Lorenzo16} the authors utilize certain gradient tracking idea to solve a constrained nonsmooth distributed problem over possibly time-varying networks.  It is worth noting that the distributed algorithms proposed in all these works converge to ss1. There has been no distributed schemes that can provably converge to ss2 for smooth non-convex problem in the form of \eqref{eq:global:consensus}.

		\section{The Gradient Primal-Dual Algorithm}
		
		In this section, we introduce the gradient primal-dual algorithm (GPDA) for solving the non-convex problem \eqref{eq:original}. Let us introduce the augmented Lagrangian (AL) as
		\begin{align}
		L(x,y) = f(x)+  \langle \lambda, Ax-b\rangle +\frac{\rho}{2}\|Ax-b\|^2\label{eq:augmented},
		\end{align}
		where $\lambda\in\mathbb{R}^{M}$ is the dual variable. The steps of the GPDA algorithm are described in the table below.
		
		Each iteration of the GPDA  performs a gradient descent step on the AL (with stepsize being $1/\beta$), followed by taking one step of approximate dual gradient ascent (with stepsize $\rho>0$). The GPDA is closely related to the classical Uzawa primal-dual method \cite{UZAWA58}, which has been utilized to solve {\it convex} saddle point problems and linearly constrained {\it convex} problems \cite{Nedic2009saddle}. It is also related to the proximal method of multipliers (Prox-MM) first developed by Rockafellar in \cite{rockafellar1976augmented}, in which a proximal term has been added to the augmented Lagrangian in order to make it strongly convex  in each iteration. The latter method has also been applied for example, in solving certain large-scale linear programs; see \cite{wright_proximal}. However the theoretical results derived for Prox-MM in \cite{rockafellar1976augmented, wright_proximal} are only developed for convex problems. Further, such an algorithm requires that the proximal Lagrangian to be optimized with {\it increasing} accuracy as the algorithm progresses. %We also note that the original prox-PDA developed in \cite{hong16decomposing} utilizes the square of a generalized $\ell_2$ norm to regularize the update step \eqref{eq:x:update}. Our analysis below also applies to that case, but for simplicity of notation we will use the version above.
		Finally, we note that both step \eqref{eq:x:update} and \eqref{eq:mu:update} can be decomposable over the variables, therefore they are easy to be implemented in a distributed manner (as will be explained shortly).
		
		\begin{center}
			\vspace{-0.5cm}
			\fbox{
				\begin{minipage}{\columnwidth}
					\smallskip
					\centerline{\bf {Algorithm 1. The gradient primal-dual algorithm}}
					\smallskip
					At iteration $0$, initialize $\lambda^0 $ and $x^0$.
					
					At each iteration $r+1$, update variables  by:
					\begin{subequations}
						\begin{align}
						x^{r+1}& =\arg\min\; \langle \nabla f(x^r) + A^T\lambda^r  + \rho A^T (A x^r-b), x-x^r\rangle +\frac{\beta}{2}\|x-x^r\|^2\label{eq:x:update}\\
						\lambda^{r+1}& = \lambda^r +\rho \left(A x^{r+1} - b \right)				\label{eq:mu:update}
						\end{align}
					\end{subequations}
					
				\end{minipage}
			}
		\end{center}

		\subsection{Application in distributed optimization problem}
		
		To see how the GPDA can be specialized to the problem of distributed optimization over the network  \eqref{eq:global:consensus:equiv}, let us begin by writing the optimality condition of  \eqref{eq:x:update}. We have
		%	\begin{subequations}
		\begin{align}\label{eq:optimality}
		\nabla f(x^r) + A^T \lambda^r + \rho A^T A x^r + \beta (x^{r+1}-x^r) = 0.
		%	\nabla f(x^{r-1}) + A^T \lambda^{r-1} + \rho A^T A x^{r-1} + \beta (x^{r}-x^{r-1}) = 0.
		\end{align}
		%	\end{subequations}	
		Subtracting \eqref{eq:optimality} with its counterpart at iteration $r$, we obtain
		\begin{align*}
		&\nabla f(x^r) - \nabla f(x^{r-1}) + A^T (\lambda^r-\lambda^{r-1}) + \rho A^T A (x^r-x^{r-1}) +\beta w^{r+1} =0.
		\end{align*}
		where we have defined $w^{r+1} = (x^{r+1}-x^r) -(x^r-x^{r-1})$.
		Rearranging, and use the fact that $A^T A =L_{-}\in\mathbb{R}^{N\times N}$ is the {\it signed Laplacian matrix}, and $b=0$ in \eqref{eq:global:consensus:equiv}, we obtain
		\begin{align}\label{eq:distributed:iteration}
		x^{r+1} & = x^r + (x^{r}-x^{r-1}) +\frac{1}{\beta}\big(-\nabla f(x^r)+ \nabla f(x^{r-1}) -\rho L_{-}x^r - \rho L_{-} (x^r-x^{r-1})\big).
		\end{align}
		Consider problem \eqref{eq:global:consensus:equiv} (for simplicity assume that $g\equiv 0$), the above iteration can be implemented in a distributed manner, where each agent $i$ performs
		\begin{align}%\label{eq:distributed:iteration2}
		x_i^{r+1} & = x_i^r + (x_i^{r}-x_i^{r-1}) +\frac{1}{\beta}\bigg(-\nabla f_i(x_i^r)+ \nabla f_i(x_i^{r-1})  -2\rho \left(d_i x^r_i-\sum_{j\in \mathcal{N}_i} x^r_j\right) + \rho \big(d_ix_i^{r-1} -\sum_{j\in\mathcal{N}_i} x_j^{r-1} \big)\bigg),\nonumber
		\end{align}
		where $\mathcal{N}_i:=\{j\mid j\ne i, (i,j)\in\mathcal{E}\}$ is the set of neighbors of node $i$; $d_i$ is the degree for node $i$. Clearly, to implement this iteration
		each node only needs to know the information from the past two iterations about its immediate neighbors. %Note that this iteration is slightly different from Eq. (24) derived in the original Prox-PDA paper \cite{hong17icml}, because here we have two parameters $\rho$ and $\beta$, while in \cite{hong17icml} only one parameter has been used.
		
		%\section{Convergence Analysis}
		\subsection{Convergence to ss1 solutions}
		We first state our main assumptions.
		\vspace{-0.2cm}
		\begin{itemize}
			\item[A1.] The function $f(x)$ is smooth and has Lipschitz continuous gradient, as well as Lipschitz continuous Hessian: %In particular, the following is true
			\begin{align}
			\hspace{-0.9cm}\|\nabla f(x) - \nabla f(y)\|&\le L\|x-y\|, \; \forall~x,y\in\mathbb{R}^N\; \label{eq:Lip}\\
			\hspace{-0.9cm}\|\nabla^2 f(x) - \nabla^2 f(y)\|&\le M \|x-y\|, \; \forall~x,y\in \mathbb{R}^N\; \label{eq:Lip:Hessian}.
			\vspace{-0.3cm}
			\end{align}
			\vspace{-0.4cm}
			\item [A2.] The function $f(x)$ is lower bounded over $x\in\mathbb{R}^N$. Without loss of generality, assume that $f(x)\ge 0$.
			\vspace{-0.2cm}
			\item [A3.] The constraint $Ax = b$ is feasible over $x\in X$.  Further, $A^T A$ is {not} full rank.
			\vspace{-0.2cm}
			\item [A4.] The function $f(x)+\frac{\rho}{2}\|Ax-b\|^2$ is coercive.
			\vspace{-0.2cm}
			\item [A5.]  The function $f$ is proper and it satisfies the Kurdyka-{\L}ojasiewicz (K{\L}) property. That is, at $\hat{x}\in\mathbb{R}$ if there exist $\eta \in (0\;\infty]$, a neighborhood $V$ of $\hat{x}$ and a continuous concave function $\phi: [0, \; \eta)\to \mathbb{R}_{+}$ such that: 1) $\phi(0) = 0$ and $\phi$ is continuously differentiable on $[0,\; \eta]$ with positive derivatives; 2) for all $x\in\mathbb{R}^N$, satisfying $f(\hat{x})<f(x)<f(\hat{x})+\eta$, it holds that
			\begin{align}
			\phi'(f(x)-f(\hat{x}))\mbox{dist}(0,\partial f(x))\ge 1.
			\end{align}
			where $\partial f(x)$ is the limiting subdifferential defined as
			\begin{align*}
			&\partial f(x) = \bigg\{\{v\in\mathbb{R}^N: \exists x^t\to x, v^t\to v, ~\mbox{with}~\lim\inf_{z\to x^t}\frac{f(x)-f(x^t)-\langle v^t,z-x^t\rangle}{\|x-x^t\|}\ge 0, \forall~t\bigg\}.
			\end{align*}
			\vspace{-0.2cm}
		\end{itemize}
		\vspace{-0.2cm}
		We comment that a wide class of functions enjoys the  K{\L} property, for example a semi-algebraic function is a KL function; for detailed discussions of the K{\L} property we refer the readers to \cite{Bolte14, Li14splitting}.
		%	We comment that Assumption [A3] makes the problem \eqref{eq:original} non-trivial. Otherwise,  if $A^T A$ is full rank, then the problem reduces to finding the unique solution of  the linear systems of equations $A x=b$.
		
		Below we will use $\sigma_{i}(\cdot)$, $\sigma_{\max}(\cdot)$, $\sigma_{\min}(\cdot)$ and $\tilde{\sigma}_{\min}(\cdot)$ to denote the $i$th, the maximum, the minimum, and the smallest {\it non-zero} eigenvalues of a matrix, respectively.
		
		The convergence of GPDA to the ss1 is similar to Theorem 3.1 in \cite{hong16decomposing} and Corollay 4.1 in \cite{hong16decomposing}. Algorithmically, the main difference is that the algorithms analyzed in \cite{hong16decomposing} do not linearize the penalty term $\frac{\rho}{2}\|A x-b\|^2$, and they make use of the same penalty and proximal parameters, that is, $\rho=\beta$.  In this work, in order to show the convergence to ss2, we need to have the freedom of tuning $\beta$ while fixing $\rho$, therefore $\beta$ and $\rho$ have to be chosen differently. However, in terms of analysis, there is no major difference between these versions. For completeness, we only outline the key proof steps in the Appendix.

		\begin{claim}\label{claim:convergence}
			Suppose Assumptions [A1] -- [A5] are satisfied. For appropriate choices of $\rho$, and $\beta$ satisfying \eqref{eq:first:order:condition} given in the appendix, and starting from any feasible point $(x^0,\lambda^0)$, the GPDA converges to the set of ss1 solutions. %satisfying
			%	\begin{align}
			%	\nabla f(x^*) +  A^T \lambda^* = 0, \quad A x^* = b.
			%	\end{align}
			
			Further, if $L(x^{r},\lambda^{r})$ is a $K{\L}$ function, then $(x^{r+1},\lambda^{r+1})$ converges globally to a unique point $(x^{*},\lambda^{*})$.
		\end{claim}
		%	
		%The following claim can be proved if needed.
		%\begin{claim}\label{claim:convergence:strong}
		%	If further $f(x)$ is semi-algebraic, then the above convergence is to a single point $(x^*,\lambda^*)$.
		%\end{claim}
		
		\subsection{Convergence to ss2}
		One can view Claim \ref{claim:convergence} as some variation of known results. On the contrary, in this section we show one of the main contributions of this work, which demonstrates that GPDA can converge to solutions beyond the ss1.  %convergence beyond converging to ss1 solutions, GPDA can
		%We begin to analyze whether the GPDA can compute the second-order stationary solutions.
		
		To this end, first let us rewrite the $x$ update step using its first-order optimality condition as follows
		\begin{align*}
		x^{r+1} = x^r - \frac{1}{\beta}\left(\nabla f(x^r)+ A^T \lambda^r + \rho A^T (A x^r-b)\right).
		\end{align*}
		Therefore the iteration can be  written as
		\begin{align*}
		&\begin{bmatrix}
		x^{r+1}\\ \lambda^{r+1}
		\end{bmatrix} = \begin{bmatrix}
		x^{r} -  \frac{1}{\beta}\left(\nabla f(x^r)+ A^T \lambda^r + \rho A^T (A x^r-b)\right)\\ \lambda^{r} + \rho (A x^{r+1}-b)
		\end{bmatrix}\nonumber\\
		&= \begin{bmatrix}
		x^{r} -  \frac{1}{\beta}\left(\nabla f(x^r)+ A^T \lambda^r + \rho A^T (A x^r-b)\right)\\
		\lambda^{r} + \rho \left(A \left(x^{r} -  \frac{1}{\beta}\left(\nabla f(x^r)+ A^T \lambda^r + \rho A^T (A x^r-b)\right)\right) - b \right).
		\end{bmatrix}
		\end{align*}
		The compact way to write the above iteration is
		\begin{align}\label{eq:recursion}
		&\begin{bmatrix}
		I_N & 0_{N\times M}\\ -\rho A & I_M
		\end{bmatrix} \begin{bmatrix}
		x^{r+1}\\ \lambda^{r+1}
		\end{bmatrix} \nonumber\\
		& = \begin{bmatrix}
		x^{r} -  \frac{1}{\beta}\left(\nabla f(x^r)+ A^T \lambda^r + \rho A^T (A x^r-b)\right)\\ \lambda^{r}-\rho b
		\end{bmatrix},
		\end{align}
		where $I_N$ denotes the $N$-by-$N$ identity matrix %{\color{red} [[[The constant term $A^T A b $ seems to be missed in the $x^r$ update rule-- The rest seems to be ok though]]]};
		$0_{N\times M}$ denotes the $N$-by-$M$ all zero matrix.
		
		Next let us consider approximating $\nabla f(x)$ near a first-order stationary solution $x^*$. Let us define
		$$H:=\nabla^2 f(x^*), \quad d^{r+1} := -x^*+x^{r+1}.$$
		Claim \ref{claim:convergence}  implies that when $\rho,\beta$ are chosen appropriately, then $d^{r+1}\to 0$. Therefore for any given $\xi>0$ there exists an iteration index $R(\xi)>0$ such that the following holds
		\begin{align}\label{eq:d:bounded}
		\|d^{r+1}\|\le \xi, \quad \forall~r-1\ge R(\xi).
		\end{align}
		
		Next let us approximate the gradients around $\nabla f(x^*)$:
		\begin{align}\label{eq:f:app}
		&\nabla f(x^{r+1}) =\nabla f(x^{*}+ d^{r+1}) \nonumber\\
		&=\nabla f(x^*) +\int^{1}_{0} \nabla^2 f(x^*+t d^{r+1})d^{r+1}  dt \nonumber\\
		& = \nabla f(x^*) +\int^{1}_{0} (\nabla^2 f(x^*+t d^{r+1}) - H)d^{r+1}  dt + H d^{r+1}\nonumber\\
		&:= \nabla f(x^*) +\Delta^{r+1}d^{r+1}   + H d^{r+1},
		\end{align}
		where in the last inequality we have defined
		\begin{align}\label{eq:Delta}
		\Delta^{r+1} := \int^{1}_{0} (\nabla^2 f(x^*+t d^{r+1}) - H) d^{r+1}dt.
		\end{align}
		From Assumption [A1] and  \eqref{eq:d:bounded} we have
		\begin{align*}
		\|\Delta^{r+1}\|\le M \|d^{r+1}\| \le M \xi, \quad \forall~r\ge R(\xi).
		\end{align*}
		Therefore we have
		\begin{align}\label{eq:delta:0}
		\lim_{r\to\infty}\|\Delta^{r+1}\|\to 0.
		\end{align}
		Using the approximation \eqref{eq:f:app}, we obtain
		\begin{align}\label{eq:diff:gradf}
		\nabla f(x^{r}) = \nabla f(x^*) +\Delta^{r}d^{r}   + H d^{r}
		\end{align}
		Plugging \eqref{eq:diff:gradf} into \eqref{eq:recursion}, the iteration \eqref{eq:recursion} can be written as
		%	{\small
		%	\begin{align}\label{eq:x:compact}
		%&	\begin{bmatrix}
		%	I_N& 0_{N\times M}\\ -\rho A & I_M
		%	\end{bmatrix} \begin{bmatrix}
		%	x^{r+1}\\ \lambda^{r+1}
		%	\end{bmatrix}  = \begin{bmatrix}
		%	I_N-  \frac{1}{\beta}\left(H + \rho A^T A \right)& -\frac{1}{\beta} A^T\\ 0_{M\times N} & I_M
		%	\end{bmatrix}\begin{bmatrix}
		%	x^{r}\\ \lambda^{r}
		%	\end{bmatrix}\nonumber\\
		%	&+ \begin{bmatrix}
		%	\nabla f(x^*) + \Delta^r d^r -H x^*\\ -\rho b
		%	\end{bmatrix}.
		%	\end{align}}
		%	Or equivalently
		{\begin{align}\label{eq:x:compact:2}
			&	\begin{bmatrix}
			x^{r+1}\\ \lambda^{r+1}
			\end{bmatrix}  = \begin{bmatrix}
			I_N & 0_{N\times M}\\ \rho A & I_M
			\end{bmatrix}\begin{bmatrix}
			I_N -  \frac{1}{\beta}\left(H + \rho A^T A \right)& -\frac{1}{\beta} A^T\\ 0_{M\times N} & I_M
			\end{bmatrix}\begin{bmatrix}
			x^{r}\\ \lambda^{r}
			\end{bmatrix}\nonumber\\
			&+\begin{bmatrix}
			I_N & 0_{N\times M}\\ \rho A & I_M
			\end{bmatrix} \begin{bmatrix}
			\nabla f(x^*) + \Delta^r d^r -H x^*\\ -\rho b
			\end{bmatrix}
			\end{align}}

		%\begin{align}\label{eq:x:compact}
		%\begin{bmatrix}
		%x^{r+1}\\ x^r
		%\end{bmatrix} &= \begin{bmatrix}
		%2 I - \frac{1}{\beta} H -\frac{2\rho}{\beta} A^T A - \frac{1}{\beta}\Delta^{r} & \quad -I +\frac{1}{\beta} H + \frac{\rho}{\beta} A^T A + \frac{1}{\beta} \Delta^{r-1} \\ I & \quad 0
		%\end{bmatrix} \begin{bmatrix}
		%x^r \\ x^{r-1}
		%\end{bmatrix} \nonumber\\
		%&\quad + \begin{bmatrix}
		%\frac{\rho}{\beta} A^T b + \frac{1}{\beta}(\Delta^r-\Delta^{r-1})x^*\\0
		%\end{bmatrix}.
		%\end{align}
		Then the above iteration can be compactly written as
		\begin{align}\label{eq:system}
		z^{r+1}  = Q^{-1}T z^r + Q^{-1}c^r
		\end{align}
		for some appropriately defined vectors $z^{r+1}, z^r, c^r$ and  matrices $M, T$ which are given below
		\begin{align}\label{eq:perturb}
		T&:=  \begin{bmatrix}
		I_N-  \frac{1}{\beta}\left(H + \rho A^T A \right)& -\frac{1}{\beta} A^T \\ 0_{M\times N} & I_M
		\end{bmatrix}\in\mathbb{R}^{(N+M)\times (N+M)}\nonumber\\
		Q &: =\begin{bmatrix}
		I_N & 0_{N\times M}\\ -\rho A & I_M
		\end{bmatrix}\in\mathbb{R}^{(N+M)\times (N+M)} \\
		c^r &:= \begin{bmatrix}
		\nabla f(x^*) + \Delta^r d^r -H x^*\\ -\rho b
		\end{bmatrix}, \; 	z := \begin{bmatrix}
		x\\ \lambda
		\end{bmatrix}
		\end{align}
		It is clear that $c^r$ is a bounded sequence.
		As a direct result of Claim \ref{claim:convergence}, we can show that every fixed point of the above iteration is an ss1 solution for problem \eqref{eq:original}.
		\begin{corollary}\label{cor:fixed:point}
			Suppose that Assumptions [A1]--[A5] are satisfied, and the parameters are chosen according to \eqref{eq:first:order:condition}. Then every fixed point of the mapping $g(z)$ defined below, is a first-order stationary solution for problem \eqref{eq:original}.
			\begin{align}\label{eq:g}
			&g(z):= g([z_1, z_2]) = \begin{bmatrix}
			I_N & 0_{N\times M}\\ \rho A & I_M
			\end{bmatrix}  \begin{bmatrix}
			z_1 -  \frac{1}{\beta}\left(\nabla f(z_1)+ A^T z_2 + \rho A^T A z_1\right)\\ z_2-\rho b
			\end{bmatrix} \nonumber. %~\mbox{with}~z_1\in\mathbb{R}^N, \; z_2\in\mathbb{R}^M.
			\end{align}
		\end{corollary}
		
		%We would like to show that around a strict saddle point \eqref{eq:strict}, the iteration expressed in \eqref{eq:x:compact} is not stable. Note that since  $c^r$ is a bounded sequence, i.e.,
		%$$\|c^r\|  = \left\|\frac{\rho}{\beta} A^T b + \frac{1}{\beta}(\Delta^r-\Delta^{r-1})x^*\right\| \le \frac{\rho}{\beta}\|A^T b\|+\frac{2M}{\beta}\|x^*\|\xi,\quad \forall~r\ge R(\xi)$$
		%then if for all $r$ large enough, matrix $T^{r+1}$ has at least one eigenvalue that lies outside of the unit circle, the system is not ``stable" around $x^*$.
		
		To proceed, we analyze the dynamics of the system \eqref{eq:system}. The following claim is a key result that characterizes the eigenvalues for the matrix $Q^{-1}T$. We refer the readers to the appendix for detailed proof.
		\begin{claim}\label{claim:sigma:T}
			Suppose Assumptions [A1] -- [A5] hold, and that
			\begin{align*}
			\beta > \sigma_{\max} (H + \rho A^T A).
			\end{align*}
			Let $(x^*, \lambda^*)$ be an ss1 solution satisfying \eqref{eq:first_order}, and that $x^*$ is a strict saddle  \eqref{eq:strict}. Let $\sigma_i(Q^{-1}T)$ be the $i$th eigenvalue for matrix $Q^{-1}T$.
			Then $Q^{-1}T$ is invertible, and there exists a real scalar $\delta^*>0$ which is independent of iteration index $r$,  such that the following holds:
			\begin{align*}
			\exists~i\in [N], \;\; \mbox{\rm s.t.}\;\; \sigma_i(Q^{-1}T) = 1+ \delta^*.
			\end{align*}
		\end{claim}

		\begin{theorem}\label{thm:main}
			Suppose that Assumptions [A1]--[A5] hold true, and that the following parameters are chosen
			\begin{align}
			\beta >\sigma_{\max}(\rho A^T A) + L, \quad \mbox{and}\quad \beta, \rho~\mbox{satisfy}~\eqref{eq:first:order:condition}.
			\end{align}	
			Suppose that $(x^{0}, \lambda^0)$ are initialized randomly. Then with probability one, the iterates $\{(x^{r+1},\lambda^{r+1})\}$ generated by the GPDA  converges  to an ss2 solution  \eqref{eq:first:and:second}.
		\end{theorem}
		\vspace{-0.3cm}
		\noindent{\bf Proof.} %The first possible way to show this is to utilize existing stability results in dynamic discrete linear systems.  {\red[?]}
		%{\red [this part requires a bit more argument.]}
		We utilize the stable manifold theorem \cite{Shub87,lee16}. We will verify the conditions given in Theorem 7  \cite{lee16} to show that the system  \eqref{eq:system} is not stable around strict saddle points.
		
		\noindent{\bf Step 1.} We will show that the mapping $g(z)$ defined in \eqref{eq:g} is diffeomorphism.

		%Finally, we investigate the question that whether the mapping $g$ from $(x^{r-1}, x^r)$ to $(x^r, x^{r+1})$ is diffeomorphism.
		First, suppose there exists $w_1 =(x_1, y_1)$, $w_2 = (x_2,y_2)$ such that $g(w_1) = g(w_2)$. Using the definition of $g$, and the fact that the matrix $[I\; 0; -\rho A \; I]$ is invertible,  we obtain $y_2 = y_1$. Using the above two results, we obtain
		\begin{align*}
		\hspace{-0.3cm}- x_1+ \frac{1}{\beta}(\rho A^T A x_1  + \nabla f (x_1)) = - x_2+ \frac{1}{\beta}(\rho A^T A x_2 + \nabla f (x_2)).
		\end{align*}
		Then we have
		\begin{align*}
		(x_1-x_2) = \frac{1}{\beta}\left(\nabla f(x_1)-\nabla f(x_2)\right) + \frac{\rho}{\beta} A^T A (x_1-x_2)
		\end{align*}
		This implies that
		\begin{align}
		\|x_1-x_2\|\le \left(\frac{L}{\beta} +\frac{\rho}{\beta}\sigma_{\max}(A^T A)\right)\|x_1-x_2\|\nonumber.
		\end{align}
		Suppose that the following is true
		\begin{align}
		\beta > \sigma_{\max}(\rho A^T A) + L.\label{eq:beta}
		\end{align}
		Then we have $x_1=x_2$, implying $y_1=y_2$. This says that the mapping $g$ is injective.
		
		To show that the mapping is surjective, we see that
		for a given tuple $(x^{r+1}, \lambda^{r+1})$, the iterate $x^{r}$ is given by
		\begin{align*}
		\ell(x^{r+1},\lambda^{r+1}) = - x^{r}+ \frac{1}{\beta}(\rho A^T A x^{r}  + \nabla f (x^{r})) %= \left(\frac{1}{\beta}\rho A^T A -I\right)x^{r-1} + \nabla f(x^{r-1})
		\end{align*}
		where $\ell(x^{r+1},\lambda^{r+1})$ is some function of $(\lambda^{r+1}, x^{r+1})$. It is clear that $x^{r}$ is the unique solution to the following convex problem [with $\beta$ satisfying \eqref{eq:beta}]
		\begin{align*}
		x^{r} = \arg\min_{x}\frac{1}{2}\|x - \ell(x^{r+1},\lambda^{r+1})\|^2 - \frac{1}{\beta}\left( f(x) + \frac{\rho}{2}\| A x\|^2\right).
		\end{align*}
		Additionally, using the definition of the mapping $g$ in \eqref{eq:g}, we have that the Jacobian matrix for the mapping $g$ is given by
		\begin{align}
		D g(z)&= \begin{bmatrix}
		I_N & 0_{N\times M}\\ -\rho A & I_M
		\end{bmatrix}\begin{bmatrix}
		I -  \frac{1}{\beta}\left(H + \rho A^T A \right)& -\frac{1}{\beta} A^T\\ 0_{M\times N} & I_M
		\end{bmatrix} \nonumber\\
		&= Q^{-1}T.
		\end{align}
		Then it has been shown in Claim \ref{claim:sigma:T} that as long as the following is true
		\begin{align}
		\beta>  {L} + {\rho}\sigma_{\max}(A^T A)
		\end{align}
		the Jacobian matrix $D g(z)$ is invertible. By applying the inverse function theorem, $g^{-1}$ is continuously differentiable.
		
		\noindent{\bf Step 2.} We can show that at a strict saddle point $x^*$, for the Jacobian matrix $Dg(z^*)$ evaluated at $z^*=(x^*,\lambda^*)$, the span of the eigenvectors corresponding to the eigenvalues of magnitude less than or equal to 1  is not the full space. This is easily done since according to Claim \ref{claim:sigma:T}, $Dg(z^*) = Q^{-1} T$ has one eigenvalue that is strictly greater than 1.
		
		\noindent{\bf Step 3.} Combining the previous two steps, and by utilizing Theorem 7 \cite{lee16}, we conclude that with random initialization, the GPDA converges to the second-order stationary solutions with probability one.  	\QED

		\section{The Gradient ADMM Algorithm}
		In this section, we extend the argument in the previous section to an algorithm belonging to the class of method called alternating direction method  of multipliers (ADMM). Although the main idea of the analysis extends those in the previous section, the presence of {\it two} blocks of primal variables instead of one significantly complicates the analysis.
		
		Consider the following problem
		\begin{align}\label{eq:original:block}
		\min\; f(x) + g(y) \quad \mbox{s.t.}\quad Ax+By = b
		\end{align}
		where $x\in\mathbb{R}^{N_1}$, $y\in\mathbb{R}^{N_2}$ and $N_1 + N_2 = N$; $b\in\mathbb{R}^{M}$. Clearly the global consensus problem \eqref{eq:global:consensus:1} can be formulated into the above two-block problem, with the following identification: $x:=\{x_1,\cdots, x_N\}$, $y:=x_0$, $f(x):=\sum_{i=1}^{N}f_i(x_i)$, $g(y):=g(x_0)$, $A=I_N$, $B=-1$, $b=0$.

		For this problem, the first- and second-order necessary conditions are given by [cf. \eqref{eq:first:and:second}]
		\begin{align}
		%	\begin{split}
		\hspace{-0.4cm}	&\nabla f(x^*)+ (\lambda^*)^T A =0,\quad \nabla g(y^*)+ (\lambda^*)^T B =0, \label{eq:second-order-admm}\\
		&z^T \begin{bmatrix}
		\nabla^2 f(x^*) & \hspace{-0.2cm}0\\
		0 & \hspace{-0.2cm}\nabla^2 g(x^*)
		\end{bmatrix} z\succeq 0,  \forall y\in \left\{z \mid \begin{bmatrix}
		A^T A & \hspace{-0.2cm}A^T B \\ B^T A & \hspace{-0.2cm}B^T B
		\end{bmatrix} z =0\right\}.\nonumber
		%	\end{split}
		\end{align}
		Similarly as before, we will refer to solutions satisfy the first line as ss1 solutions, and those that satisfy both as ss2 solutions.
		Therefore, a strict saddle point is defined as a point $(x^*,y^*,\lambda^*)$ that satisfies the following conditions
		\begin{align}
		%	\begin{split}
		&\nabla f(x^*)+ (\lambda^*)^T A =0,\quad \nabla g(y^*)+ (\lambda^*)^T B =0,\nonumber\\
		&z^T \begin{bmatrix}
		\nabla^2 f(x^*) & 0\\
		0 & \nabla^2 g(y^*)
		\end{bmatrix} z\le -\sigma\|z\|^2,\; \mbox{for some}\; \sigma>0, \; z \; \mbox{satisfying} \; \begin{bmatrix}
		A^T A & A^T B \\ B^T A & B^T B
		\end{bmatrix} z = 0.\label{eq:strict:admm}
		%	\end{split}
		\end{align}
		Define the AL function as
		\begin{align*}
		\hspace{-0.3cm}L(x,y;\lambda) = f(x)+g(y) +\langle \lambda, Ax+By-b\rangle +\frac{\rho}{2}\|Ax +By-b\|^2.
		\end{align*}
		The gradient ADMM  (G-ADMM) algorithm that we propose is given below.
		\begin{center}
			\vspace{-0.2cm}
			\fbox{
				\begin{minipage}{\columnwidth}
					\smallskip
					\centerline{\bf {Algorithm 2. The gradient ADMM}}
					\smallskip
					At iteration $0$, initialize $\lambda^0 $ and $x^0$.
					
					At each iteration $r+1$, update variables  by:
					\begin{subequations}
						\begin{align}
						x^{r+1}& =\arg\min_{x}\; \langle \nabla f(x^r) + A^T \lambda^r + \rho A^T (A x^r + By^r-b),  x-x^r\rangle +\frac{\beta}{2}\|x-x^r\|^2\label{eq:x:update:admm}\\
						y^{r+1}& =\arg\min_{y}\; \langle \nabla g(y^r) + B^T \lambda^r + \rho B^T (A x^{r+1} + By^r-b), y-y^r\rangle +\frac{\beta}{2}\|y-y^r\|^2 \label{eq:y:update:admm}\\
						\lambda^{r+1}& = \lambda^r +\rho \left(A x^{r+1} + By^{r+1} - b \right)				\label{eq:mu:update:admm}.
						\end{align}
					\end{subequations}
					
				\end{minipage}
			}
		\end{center}
		We note that in the GADMM, the $x$ and $y$ steps perform gradient steps to optimize the AL, instead of performing the exact minimization as the original convex version of ADMM does \cite{BoydADMMsurvey2011,EcksteinBertsekas1992}. The reason is that the direct minimization may not be possible because the non-convexity of $f$ and $g$ makes the subproblem of minimizing the AL w.r.t. $x$ and $y$ also non-convex. Note that the gradient steps have been used  in the primal updates of ADMM when dealing with convex problems,  see \cite{gao14}, but their analyses do not extend to the non-convex setting. %However to the best of our knowledge there has been no non-convex ADMM that utilizes
		
		It is also worth noting that the key difference between Algorithm 2 and 1 is that, in the $y$ update step \eqref{eq:y:update:admm} of Algorithm 2, the newly updated $x^{r+1}$ is used. If in this step $x^{r}$ is used instead of $x^{r+1}$, then Algorithm 2 is equivalent to Algorithm 1. Also there are quite a few recent works applying ADMM-type method to solve a number of non-convex problems; see, e.g., \cite{li14nonconvex, wang15nonconvexadmm, Goncalves17} and the references therein. However, to the best of our knowledge, these algorithms do not take  exactly the same form as Algorithm 2 described above, despite the fact that their analyses  all appear to be quite similar (i.e., some potential function based on the AL is shown to be descending at each iteration of the algorithm). In particular, in \cite{Goncalves17}, both the $x$ and $y$ subproblems are solved using a proximal point method; In \cite{jiang16admm}, the $x$-step is solved using the gradient step, while the $y$-step is solved using the conventional exact minimization. Of course, none of these works analyzed the convergence of these methods to ss2 solutions.
		
		\subsection{Application in global consensus problem}
		We discuss how Algorithm 2 can be applied to solve the global consensus \eqref{eq:global:consensus:1}. For this problem, the distributed nodes and the master node alternate between their updates: 
		\begin{subequations}
			\begin{align*}
			x^{r+1}_i & = \arg\min_{x_i} \; \left\langle \nabla f_i(x^r_i) + \lambda^r_i + \rho (x^r_i - x^r_0), x_i-x^r_i \right\rangle +\frac{\beta}{2}\|x_i-x^r_i\|^2, \; \forall~i\\
			x^{r+1}_0 & = \arg\min_{x_0}\; \langle \nabla g(x_0)  - \sum_{i=1}^{N}(\lambda^r_i + \rho (x^{r+1}_i - x^r_0)), x_0-x^r_0\rangle +\frac{\beta}{2}\|x_0-x^r_0\|^2.
			\end{align*}
		\end{subequations}
		Clearly, for fixed $x_0$, the distributed nodes are able to perform their computation completely in parallel. %These iterations are closely related to the deterministic version of Algorithm 3 in \cite{hong14nonconvex_admm}, with the difference that in our current scheme the penalty term $\rho\|x_0-x_i\|^2$ has also been linearized. The algorithms proposed in \cite{hong14nonconvex_admm} represent the first ones that provide global sublinear rate for non-convex global consensus problem (to first-order stationary solutions), but it is not clear whether similar schemes can achieve second-order stationary solutions as well.

		\subsection{Convergence to first-order stationary solutions}
		First we make the following assumptions.
		\begin{itemize}
			\vspace{-0.3cm}
			\item[B1.] The function $f(x)$ and $g(y)$ are smooth and both have Lipschitz continuous gradient and Hessian, with constants $L_f$, $L_g$, $M_f$ and $M_g$.
			%	\begin{align}
			%	&\|\nabla f(x) - \nabla f(y)\|\le L_f\|x-y\|, \quad \forall~x,y\in\mathbb{R}^n\; \label{eq:Lip}\\
			%	&\|\nabla^2 f(x) - \nabla^2 f(y)\|\le M_f \|x-y\|, \quad \forall~x,y\in \mathbb{R}^n\; \label{eq:Lip:Hessian}.
			%	\end{align}
			\vspace{-0.2cm}
			\item [B2.] $f(x)$ and $g(y)$ are lower bounded over $\mathbb{R}^N$. Without loss of generality, assume $f(x)\ge 0, g(y)\ge 0$.
			\vspace{-0.2cm}
			\item [B3.] $Ax +By = b$ is feasible over $x\in \dom(f)$ and $y\in \dom(g)$; the matrix $[A; B]\in\mathbb{R}^{M\times N}$ is {\it not} full rank.
			\vspace{-0.2cm}
			\item [B4.] $f(x) + g(y)+\frac{\rho}{2}\|Ax +By-b\|^2$ is a coercive function.
			\vspace{-0.2cm}
			\item [B5.] $f(x)+g(x)$ is a (K{\L}) function given in [A5].
		\end{itemize}
		\vspace{-0.3cm}
		Based on the above assumptions, the convergence of Algorithm 2 to the ss1 solutions can be shown following similar line of arguments as in \cite{li14nonconvex, wang15nonconvexadmm, Goncalves17,jiang16admm}. However, since the exact form of this algorithm has not appeared before, for completeness we provide the proof outline in the appendix.
		\begin{claim}\label{claim:convergence:admm}
			Suppose Assumptions [B1] -- [B5] are satisfied. For appropriate choices of $\beta, \rho$ [see \eqref{eq:condition:admm} in the Appendix for the precise expression], and starting from any point $(x^0, y^0, \lambda^0)$, Algorithm 2 converges to the set of ss1 points. %satisfying{\small
			%\begin{align*}
			%\nabla f(x^*) +  A^T \lambda^* = 0,  \nabla g(y^*) + B^T \lambda^* = 0,  A x^* + By^*= b.
			%\end{align*}}
			Further, if $L(x^{r+1}, y^{r+1},\lambda^{r+1})$ is a K{\L} function, then Algorithm 2 converges globally to a unique point $(x^{*}, y^*, \lambda^{*})$.
		\end{claim}
		
		\subsection{Convergence to ss2 solutions}
		The optimality conditions for the $(x,y)$ update is given as
		\begin{align*}
		&\nabla f(x^r) +  A^T \lambda^r + \rho A^T (A x^r+ B y^r-b) + \beta (x^{r+1}-x^r)=0\nonumber\\
		&\nabla g(y^r) + B^T \lambda^r + \rho B^T (A x^{r+1}+ B y^r-b) + \beta (y^{r+1}-y^r)=0\nonumber.
		\end{align*}
		These conditions combined with the update rule of the dual variable give the following compact form of the algorithm
		\begin{align*}%\label{eq:admm:iterations}
		\hspace{-0.5cm}\begin{bmatrix}
		x^{r+1}\\y^{r+1}\\\lambda^{r+1}
		\end{bmatrix} = \begin{bmatrix}
		x^r-\frac{1}{\beta}\left(\nabla f(x^r) + A^T \lambda^r + \rho A^T (A x^r+ By^r -b)\right)\\
		y^r-\frac{1}{\beta}\left(\nabla g(y^r) + B^T \lambda^r + \rho B^T (A x^{r+1}+ By^r -b)\right)\\
		\lambda^r + \rho\left(A x^{r+1} + B y^{r+1} -b\right)
		\end{bmatrix}.
		\end{align*}
		%These conditions imply that
		%\begin{align*}
		%&\nabla f(x^r)-\nabla f(x^{r-1}) + A^T (\lambda^r-\lambda^{r-1}) \nonumber\\
		%&\quad + \rho A^T \left(A(x^r-x^{r-1})+ B(y^r-y^{r-1})\right) + \beta [(x^{r+1}-x^r)-(x^r-x^{r-1})]=0\nonumber\\
		%&\nabla g(y^r)-\nabla g(y^{r-1}) + B^T (\lambda^r-\lambda^{r-1}) \nonumber\\
		%&\quad + \rho B^T \left(A(x^{r+1}-x^{r})+ B(y^r-y^{r-1})\right) + \beta [(y^{r+1}-y^r)-(y^r-y^{r-1})]=0\nonumber.
		%\end{align*}
		%Therefore, the new iterate $(x^{r+1}, y^{r+1})$ is given by
		%\begin{subequations}\label{eq:admm:iterations}
		%\begin{align}
		%x^{r+1} &= x^r+ (x^r-x^{r-1}) -\frac{1}{\beta}\bigg[ \nabla f(x^r) - \nabla f(x^{r-1}) + \rho A^T (A x^r + By^r -b)\nonumber\\
		%&\quad + \rho A^T\left( A(x^r-x^{r-1})+B(y^r-y^{r-1})\right)\bigg]\label{eq:x:iter}\\
		%y^{r+1} &= y^r+ (y^r-y^{r-1}) -\frac{1}{\beta}\bigg[ \nabla g(y^r) - \nabla g(y^{r-1}) + \rho B^T (A x^r + By^r -b)\nonumber\\
		%&\quad + \rho B^T\left( A(x^{r+1}-x^{r})+B(y^r-y^{r-1})\right)\bigg]\label{eq:y:iter}.
		%\end{align}
		%\end{subequations}
		To compactly write the iterations in the form of a linear dynamic system, define
		\begin{align*}
		z^{r+1}: = [x^{r+1}; y^{r+1}; \lambda^{r+1}]\in\mathbb{R}^{2N+M}.
		\end{align*}
		
		Next we approximate the iteration around a stationary solution $x^*$. Suppose that $\nabla^2 f(x^*) = H$ and $\nabla^2 g(y^*) = G$. Then similarly as the derivation of \eqref{eq:x:compact:2}, we can write
		\begin{align*}
		P z^{r+1} = T^{r} z^r + d = (T+E^{r}) z^r + d^r
		\end{align*}
		where we have defined 
		\begin{subequations}
			\begin{align}
			\hspace{-0.5cm}P &:= \begin{bmatrix}
			I_N & 0 & 0 \\\frac{\rho}{\beta} B^T A & I_N & 0 \\ -\rho A & -\rho B & I_M
			\end{bmatrix}, \; 	E^r := \begin{bmatrix}
			\Delta^{r}_{H}\\\Delta^{r}_G\\0
			\end{bmatrix}\\%\in\mathbb{R}^{(2N+M)\times (2N+M)},\\
			\hspace{-0.5cm} d &:= \begin{bmatrix}
			\frac{\rho}{\beta} A^T b +\nabla f(x^*) - \Delta^r_H x^* - H x^*\\
			\frac{\rho}{\beta} B^T b +\nabla g(y^*) - \Delta^r_G x^*- G y^*\\ -\rho b
			\end{bmatrix} \label{eq:M}\\%\in\mathbb{R}^{2N+M}\label{eq:M}\\
			\hspace{-0.5cm}T&:=\begin{bmatrix}
			I_N -\frac{1}{\beta}H-\frac{\rho}{\beta}A^T A & \hspace{-0.5cm}-\frac{\rho}{\beta} A^T B & \hspace{-0.3cm} -\frac{1}{\beta} A^T \\
			0 & \hspace{-0.5cm} I_N-\frac{1}{\beta} G +\frac{\rho}{\beta}B^T B & \hspace{-0.2cm}-\frac{1}{\beta}B^T\\
			0 & \hspace{-0.3cm} 0 & \hspace{-0.2cm} I_M
			\end{bmatrix}\label{eq:T}%\in\mathbb{R}^{(2N+M)\times(2N+M)}\label{eq:T}\\
			%E^r & = \frac{1}{\beta}\begin{bmatrix}
			%-\Delta^r_{H} & 0 &  \Delta^{r-1}_H & 0 \\ 0 & - \Delta^r_G  & 0 & \Delta^{r-1}_G\\
			%0 & 0 & 0 & 0\\
			%0 & 0 & 0 & 0
			%\end{bmatrix}, \quad
			%d^r := \begin{bmatrix}
			%\frac{\rho}{\beta} A^T b  +\frac{1}{\beta} (\Delta^r_H -\Delta^{r-1}_H) x^*\\ \frac{\rho}{\beta} B^T b +\frac{1}{\beta} (\Delta^r_G -\Delta^{r-1}_G) y^*\\ 0 \\ 0
			%\end{bmatrix}.
			\end{align}
		\end{subequations}
		with the following 
		\begin{align*}
		&\Delta_H^{r+1} := \int^{1}_{0} (\nabla^2 f(x^*+t d_x^{r+1}) - H) d^{r+1}_x dt\\
		&\Delta_G^{r+1} := \int^{1}_{0} (\nabla^2 g(y^*+t d_y^{r+1}) - G) d^{r+1}_y dt,\\
		&\mbox{with}\quad d_{x}^{r+1} := -x^*+x^{r+1},  \quad d_{y}^{r+1} := -y^*+y^{r+1}.
		\end{align*}
		By noting that $P$ is an invertible matrix,
		we conclude that the new iteration $z^{r+1}$ can be expressed as
		\begin{align}
		z^{r+1} = P^{-1}(T+ E^{r+1}) z^r + P^{-1}d^r\label{eq:z:iteration}.
		\end{align}
		
		%\begin{corollary}
		%	Suppose that Assumption [B1] -- [B3] are satisfied. Define the following mapping
		%	\begin{align}
		%	h(z) = h(z_1,z_2,z_3,z_4) = \begin{bmatrix}
		%	I & 0 & 0 & 0 \\-\frac{\rho}{\beta} B^T A & I & 0 & 0\\ 0 & 0 & I & 0\\ 0 & 0& 0 & I
		%	\end{bmatrix} \begin{bmatrix}
		%	2 I -\frac{\beta}{\rho}& & & \\ & & & \\ 0 & 0 & I & 0\\ 0 & 0 & 0 & I
		%	\end{bmatrix}
		%	\end{align}
		%\end{corollary}
		
		Now in order to analyze the stability at a point $(x^*,y^*)$, similarly as before we need to analyze the eigenvalues of the matrix $P^{-1}T$ at a stationary solution.
		
		We note that $P$ is a lower triangular matrix and $\det P=1$. This implies that $\det(P^{-1} T-\mu I) =\det(T-\mu P)$. We have the following characterization on the determinant of $T-\mu P$;  please see Appendix for detailed proof.	
		\begin{claim}\label{claim:real}
			We have the following for $\det[T-\mu P]$:	
				
			\noindent	1) $\det[T-P]=0$, i.e., $1$ is an eigenvalue of $P^{-1}T$.
			
			\noindent	2) Suppose that the following condition is satisfied
			\begin{align*}
			\beta> \rho\sigma_{\max}(A^T A) + L_f, \quad \beta> \rho\sigma_{\max}(B^T B) + L_g,
			\end{align*}
			Then $\det[T]\ne 0$, i.e., the matrix $P^{-1}T$ is invertible.
			
			\noindent	3) Define a $2N\times 2N$ matrix $U(\mu) = [U_{11}(\mu)\; U_{12}(\mu); U_{12}(\mu)\; U_{22}(\mu)]$, with
			\begin{subequations}
				\begin{align}\label{eq:Q}
				U_{11}(\mu) &= -\mu\left(2I - \frac{2\rho}{\beta}A^T A-\frac{1}{\beta} H -\mu I\right) + I - \frac{\rho}{\beta}A^T A -\frac{1}{\beta} H\\
				U_{12}(\mu) & = \mu \frac{2\rho}{\beta}A^T B -\frac{\rho}{\beta}A^T B = (2\mu-1)\frac{\rho}{\beta}A^T B \\
				U_{21}(\mu) & = \mu^2\frac{\rho}{\beta}B^T A\\
				U_{22}(\mu) & = -\mu\left(2 I -\frac{1}{\beta}G -\frac{2\rho}{\beta}B^T B -\mu I\right) + I-\frac{1}{\beta} G- \frac{\rho}{\beta}B^T B.
				\end{align}
			\end{subequations}
			Then we have $\det[U(\mu)] = \det[T-\mu P]$, and that for any $\delta\in\mathbb{R}_{+}$  the eigenvalues of $U(1+\delta)$ are the same as those of the following symmetric matrix
			\begin{align}\label{eq:symmetric}
			\begin{bmatrix}
			U_{11}(1+\delta) & \hspace{-0.3cm}{(\delta+1)}{\sqrt{2\delta+1}}\frac{\rho}{\beta}A^T B\\ {(\delta+1)}{\sqrt{2\delta+1}} \frac{\rho}{\beta}B^T A & \hspace{-0.3cm} U_{22}(1+\delta).
			\end{bmatrix}
			\end{align}
			%	where $U_{11}, U_{22}$ are two submatrices of $U$ defined in \eqref{eq:Q}.
			%	Further, for any $\delta>0$ small enough, $F(\delta)$ has at least one negative eigenvalue.
		\end{claim}

		%Let us then consider the case where $|\mu|>1$. In fact,
		Based on Claim \ref{claim:real}, we will show that the matrix $P^{-1}T$ has a {\it real} eigenvalue $\mu=1+\delta$, with $\delta>0$ being a positive number. To this end, plugging $\mu=1+\delta$ to the expression of the $U$ matrix in \eqref{eq:Q} we have
		\begin{align*}
		U_{11}(1+\delta) &= \delta^2 I + \frac{\rho}{\beta}(1+2\delta) A^T A + \frac{\delta}{\beta} H\\
		U_{21}(1+\delta) & = (1+\delta)^2 \frac{\rho}{\beta} B^T A , \quad U_{12}(1+\delta)  = (1+2\delta) \frac{\rho}{\beta}A^T B\\
		U_{22}(1+\delta) & = \delta^2 I + \frac{\rho}{\beta}(1+2\delta) B^T B + \frac{\delta}{\beta} G.
		\end{align*}
		Therefore, in this case we can express $U(1+\delta)$ as
		\begin{align*}
		U(1+\delta) = (2\delta +1) U(1) + \frac{\delta}{\beta} \begin{bmatrix}
		H & 0\\ 0 & G
		\end{bmatrix} + \delta^2\begin{bmatrix}
		I & 0 \\ \frac{\rho}{\beta}B^T A & I
		\end{bmatrix} \nonumber.
		\end{align*}
		It remains to show that there exists $\delta^*>0$ such that the determinant of the above matrix is zero. To this end, we rewrite the above expression as follows
		\begin{align}\label{eq:Q:delta}
		U(1+\delta) &= \delta\left( \frac{2\delta+1}{\delta} U(1) + \frac{1}{\beta}   \begin{bmatrix}
		H & 0\\ 0 & G
		\end{bmatrix} + \delta \begin{bmatrix}
		I & 0 \\ \frac{\rho}{\beta}B^T A & I
		\end{bmatrix}\right)\nonumber\\
		&:= \delta\left(F(\delta) + E(\delta)\right)
		\end{align}
		where for notational simplicity, we have defined 
		\begin{align*}%\label{eq:F}
		F(\delta)& = \frac{(2\delta +1)}{\delta} U(1) + \frac{1}{\beta} \begin{bmatrix}
		H & 0\\ 0 & G
		\end{bmatrix} , \;E(\delta) =  \delta\begin{bmatrix}
		I & 0 \\ \frac{\rho}{\beta}B^T A & I
		\end{bmatrix}. %\quad \tilde{A}(\delta) = A(\delta)+E(\delta).
		\end{align*}
		Note that  from \eqref{eq:strict:admm}, we know that at a strict saddle point, there exists $y$ such that
		\begin{align}
		U(1) y =0, \quad  y^T\begin{bmatrix}
		H & 0\\ 0 & G
		\end{bmatrix} y\le -\sigma\|y\|^2,
		\end{align}
		which implies
		\begin{align}
		y^T \left(\gamma U(1) + \begin{bmatrix}
		H & 0\\ 0 & G
		\end{bmatrix} \right)y \le -\sigma\|y\|^2, \; \forall~\gamma.
		\end{align}
		This further implies that the matrix  $F(\delta)$ has eigenvalue no greater than $-\sigma/\beta$ for any $\delta$.
		
		Next we invoke a matrix perturbation result \cite{Stewart90} to argue that the matrix $F(\delta)+ E(\delta)$ also has negative eigenvalue as long as the parameter $\delta>0$ is small enough.
		
		For a given matrix $\tilde{F} = F+ E\in\mathbb{R}^{N\times N}$,  let us define the following quantity, which is referred to as the optimal matching distance between $F$ and $\tilde{F}$ [see Chapter 4, Section 1, Definition 1.2 in \cite{Stewart90}]
		\begin{align}\label{eq:md}
		\mbox{md}(F, \tilde{F}): = \min_{\Pi} \max_{j\in[N]}| \tilde{\sigma}_{\Pi(j)}-\sigma_j|
		\end{align}
		where $\Pi$ is taken over all permutations of $[N]$, and $\sigma_j$ (resp $\tilde{\sigma}_j$) is the $j$th eigenvalue of $F$ (resp. $\tilde{F}$).
		We have the following  results characterizing the matching distance of two matrices $F$ and $\tilde{F}$ \cite{Stewart90}:
		\begin{claim}\label{claim:md1}
			Suppose that $F$ is diagonalizable, i.e., $X^{-1} F X =\Upsilon$. Then the following is true
			\begin{align}
			\mbox{\rm md}(F, \tilde{F})\le (2N-1)\|X\|\|X^{-1}\|\|E\|.
			\end{align}
		\end{claim}
		%  We have the following bounds on the optimal matching distance between $F$ and $\tilde{F}$ (which are not necessarily diagonalizable) \cite[Chapter 4, Section 3.2, Theorem 3.3]{Stewart90}.
		%  \begin{claim}\label{claim:md2}
		%  	The following is true
		%  	\begin{align}
		%  	\mbox{\rm md}(F, \tilde{F})\le (2N-1) (\|F\|_2 + \|\tilde{F}\|_2\|)^{1-\frac{1}{N}}\|E\|^{\frac{1}{N}}_2.
		%  	\end{align}
		%  \end{claim}
		Let us apply Claim \ref{claim:md1} to the matrices $F(\delta)$ and $F(\delta)+E(\delta)$. Note that
		\begin{align*}
		\|E\|_2 %&= \sqrt{\lambda_{\max}(E^T E)} \nonumber\\
		&= \delta\sigma_{\max}\left(\begin{bmatrix}
		I & \frac{\rho}{\beta}A^T B\\ \frac{\rho}{\beta}B^T A & \frac{\rho^2}{\beta^2} B^T A A^T B +I
		\end{bmatrix}\right) := \delta d
		\end{align*}
		where $d$ is a fixed number independent of $\delta$.
		By applying Claim \ref{claim:md1}, and using the fact that $\|X\|=1$, we obtain the following
		\begin{align}
		\mbox{md}(F(\delta),F(\delta)+E(\delta))\le (2N-1) \delta d.
		\end{align}
		%{\color{red}[[[What happens to $\|X\|$?]]]}
		
		Clearly, we can pick $\delta = \frac{\sigma}{2 d \beta(2N-1)} $, which implies that  %by utilizing the fact that $\Delta^r\to 0$, we can always find an index $R(\delta^*)>R(c)$ such that
		\begin{align}
		\mbox{md}(F(\delta), F(\delta)+ E(\delta))\le \frac{\sigma}{2\beta}.
		\end{align}
		This combined with the fact that  $F(\delta)$ has an eigenvalue smaller or equal to $-\sigma/\beta$ regardless of the choice of $\delta$, and that all the eigenvalues of $F(\delta)+ E(\delta)$ are real (cf. Claim \ref{claim:real}), we conclude  that there exists an index $i\in [N]$ such that
		\begin{align}
		\sigma_{i}(F(\delta)+ E(\delta))\le -\frac{\sigma}{2\beta}.
		\end{align}
		This implies that
			\begin{align*}
			\sigma_i(U(1+\delta)) &\stackrel{\eqref{eq:Q:delta}}= \delta\sigma_{i}(F(\delta)+ E(\delta)) \le -\frac{\sigma\delta}{2\beta} = -\frac{\sigma^2}{4 \beta^2(2N-1)}.
			\end{align*}
		In conclusion, we have the following claim.
		\begin{claim}\label{eq:positive:negative}
			There exists $\hat{\delta}>0$ and $\tilde\delta>0$ such that
			\begin{align}
			\sigma_{\min}(U(1+\hat{\delta}))<0,\;	\sigma_i (U(1+\tilde{\delta}))>1, \quad \forall~i.
			\end{align}
		\end{claim}
		\vspace{-0.2cm}
		\noindent{\bf Proof.}  The first claim comes directly form our above discussion. The second claim is also easy to see by analyzing the eigenvalues for the symmetric matrix in \eqref{eq:symmetric}, for large positive $\delta$.
		\QED
		\vspace{-0.2cm}
		Using the results in Claim  \ref{claim:real} and Claim \ref{eq:positive:negative}, and using the fact that the eigenvalues for $U(1+\delta)$ are continuous functions of $\delta$, we conclude that there exists $\delta^*>0$ such that $\det[U(1+\delta^*)]=0$.  The result below summarizes the proceeding discussion.
		\begin{claim}\label{claim:sigma:MT}
			Suppose Assumptions [B1] --[B5] hold true. Let $(x^*, y^*, \lambda^*)$ be a first-order stationary solution satisfying \eqref{eq:first_order}, and that it is a strict saddle point satisfying \eqref{eq:strict:admm}. Let $\sigma_i(P^{-1}T)$ be the $i$th eigenvalue for matrix $P^{-1}T$.
			Then the following holds:
			\begin{align}
			\exists~i\in [N], \;\; \mbox{\rm s.t.}\;\; |\sigma_i(P^{-1}T)|>1.
			\end{align}
			Further, when $\beta$ satisfies
			\begin{align} \label{eq:beta:admm:invertible}
			\hspace{-0.2cm}\beta> \rho\sigma_{\max}(A^T A) + L_f, \; \beta> \rho\sigma_{\max}(B^T B) + L_g,
			\end{align}
			The matrix $P^{-1}T$ is invertible.
		\end{claim}

		The rest of the proof uses a similar argument as  in Theorem \ref{thm:main}. We have the following result for the GADMM algorithm.
		
		\begin{theorem}\label{thm:main:2}
			Suppose that Assumptions [B1] -- [B5] hold, and $\beta, \rho$ are chosen according to \eqref{eq:beta:admm:invertible} and \eqref{eq:condition:admm} in the Appendix.
			Suppose that $(x^{0}, y^{0}, \lambda^0)$ are initialized randomly. Then with probability one, the iterates  generated by the GADMM converge  to an ss2 solution satisfying \eqref{eq:second-order-admm}.
		\end{theorem}

		\section{Conclusion}
		The main contribution of this work is to show that primal-dual based first-order methods are capable of converging to second-order stationary solutions, for linearly constrained non-convex problems. The main techniques that we have leveraged is the Stable Manifold Theorem and its recently developed connection to first-order optimization methods.  One important implication of our result is that, properly designed distributed non-convex optimization methods (for both the global consensus problem and the distributed optimization problem over a multi-agent network) can also converge to second-order stationary solutions.  To the best of our knowledge, this is the first algorithm for non-convex distributed optimization that is capable of computing second-order stationary solutions. Some preliminary numerical results (included in the appendix) also show that the proposed algorithms work well and they are able to avoid strict saddle points. 

		\appendix
		
	%	\newpage
	%	\onecolumn
		
		\section{Appendix}

		\subsection{Proof outline for Claim \ref{claim:convergence}}
		In this subsection we outline the proof steps.
		
		\noindent{\bf Step 1.} Let us define
		\begin{align}
		C_1 :=\beta I -\rho A^T A \succ 0.
		\end{align}
		By utilizing \eqref{eq:optimality} and the $\lambda$ update rule \eqref{eq:mu:update} we have
		\begin{align*}
		\nabla f(x^r) + A^T\lambda^{r+1} + \rho A^T A (x^r-x^{r+1}) + \beta (x^{r+1}-x^r)=0.
		\end{align*}
		Subtracting the above equality with the same one from the previous iteration, we obtain
		\begin{align*}
		&\nabla f(x^r) -\nabla f(x^{r+1}) + A^T(\lambda^{r+1}-\lambda^r)  + \rho A^T A ((x^r-x^{r+1})-(x^{r-1}-x^r) )+ \beta ((x^{r+1}-x^r) -(x^r-x^{r-1}))=0.
		\end{align*}
		Utilizing the fact that $\lambda^{r+1}-\lambda^r$ lies in the column space of $A$,  it is easy to show that the following inequality is true
		\begin{align}\label{eq:mu:difference:bound:1}
		\frac{1}{\rho}\|\lambda^{r+1}-\lambda^{r}\|^2&\le \frac{2L^2}{\rho\tilde{\sigma}_{\min}(A^T A)}\left\| x^{r}-x^{r-1}\right\|^2+\frac{2}{\rho\tilde{\sigma}_{\min}(A^T A)}\left\|(x^{r+1}-x^r) -(x^{r}-x^{r-1})\right\|_{C^T_1 C_1}^2\nonumber\\
		&:= \frac{2L^2}{\rho\tilde{\sigma}_{\min}(A^T A)}\left\| x^{r}-x^{r-1}\right\|^2+\frac{2}{\rho\tilde{\sigma}_{\min}(A^T A)}\left\|w^{r+1}\right\|_{C^T_1 C_1}^2
		\end{align}
		where we have defined $w^{r+1}:=(x^{r+1}-x^r) -(x^{r}-x^{r-1})$.
		%where $\theta>0$ is a constant that satisfies the following relation
		%\begin{align}
		%\theta=\sigma_{\max}\left(C_1\right).
		%\end{align}
		%Note that such a positive constant exists for sufficiently large $\beta$.  For example we can choose $\beta$ large enough and choose $\theta>0$ as
		%\begin{align}\label{eq:theta}
		%\theta = \beta- 2\rho\sigma_{\max}(A^T A).
		%\end{align}
		
		\noindent{\bf Step 2.} We have the following optimality condition for the $x$-update step:
		\begin{subequations}\label{eq:optimality:2}
			\begin{align}
			\langle \nabla f(x^r) + A^T \lambda^r + \rho A^T A x^r + \beta (x^{r+1}-x^r), x-x^{r+1}\rangle  \ge 0,\; \forall~x\\
			\langle \nabla f(x^{r-1}) + A^T \lambda^{r-1} + \rho A^T A x^{r-1} + \beta (x^{r}-x^{r-1}), x-x^r\rangle  \ge 0, \; \forall~x.
			\end{align}
		\end{subequations}
		Plugging $x=x^{r}$ into the first inequality and $x=x^{r+1}$ into the second, and subtracting the two inequalities, we obtain
		\begin{align}
		\rho\langle A^T (A x^{r+1} -b),  x^{r+1}-x^r\rangle \le \langle \nabla f(x^{r-1})- \nabla f(x^r) - C_1w^{r+1},x^{r+1}-x^r\rangle.
		\end{align}
		This implies that
		\begin{align}
		&\frac{\rho}{2} \|A x^{r+1}-b\|^2+\frac{1}{2}\|x^{r+1}-x^r\|^2_{C_1}\nonumber\\
		&\le \frac{\rho}{2} \|A x^{r}-b\|^2+\frac{1}{2}\|x^{r}-x^{r-1}\|^2_{C_1} + \frac{L}{2}\|x^{r+1}-x^r\|^2+ \frac{L}{2}\|x^r-x^{r-1}\|^2 -\frac{1}{2}\|w^{r+1}\|^2_{C_1}.
		\end{align}
		
		Second, we can show that after one primal dual step the AL descends in the following manner
		\begin{align}
		L(x^{r+1},\lambda^{r+1}) - L(x^r,x^{r})\le -\frac{\beta}{2}\|x^{r+1}-x^r\|^2+\frac{1}{\rho}\|\lambda^{r+1}-\lambda^r\|^2
		\end{align}
		whenever $\beta$ is chosen to satisfy
		\begin{align}
		\beta>L+\sigma_{\max}(\rho A^T A).
		\end{align}
		Therefore, combining the previous two inequalities, and using the result in Step 1,  we can show that the following inequality is true (for some constant $c>0$ to be chosen later)
		\begin{align}\label{eq:lag:bound:1}
		&L(x^{r+1},\lambda^{r+1}) + \frac{c\rho}{2}\|Ax^{r+1}-b\|^2+\frac{c}{2}\|x^{r+1}-x^r\|_{C_1}^2+\left(\frac{2L^2}{\rho\tilde{\sigma}_{\min}(A^T A)}+\frac{cL}{2}\right)\|x^{r+1}-x^r\|^2 \\
		&\le L(x^r, \lambda^r)+ \frac{c\rho}{2}\|Ax^{r}-b\|^2+ \frac{c}{2}\|x^{r}-x^{r-1}\|_{C_1}^2 +\left(\frac{2L^2}{\rho\tilde{\sigma}_{\min}(A^T A)}+\frac{cL}{2}\right)\|x^{r-1}-x^r\|^2\nonumber\\
		&\quad -\left(\frac{\beta}{2}- {cL}- \frac{2L^2}{\rho\tilde{\sigma}_{\min}(A^T A)}\right)\|x^{r+1}-x^r\|^2  - (w^{r+1})^T \left({\frac{c}{2} C_1 - \frac{2 C_1^T C_1}{\rho\tilde{\sigma}_{\min}(A^T A)}}\right) w^{r+1} \nonumber
		\end{align}
		where $\tilde{\sigma}_{\min}(A^T A)$ is the smallest non-zero eigenvalue for $A^T A$.

		\noindent{\bf Step 3.} It is easy to show that $L(x^{r+1},\lambda^{r+1})$ is lower bounded; see Lemma 3.5 in \cite{hong16decomposing}.
		
		\noindent{\bf Step 4.} It is also easy to show that there exists a constant $c(\beta,\rho,\sigma_{\max})>0$ (which is a function of $\beta, \rho, \sigma_{\max}$) such that the following holds
		\begin{align}
		\|\nabla L(x^{r+1},\lambda^{r+1})\|\le c(\beta, \rho,\sigma_{\max})(\|x^{r+1}-x^r\| + \|x^{r}-x^{r-1}\|).
		\end{align}
		
		\noindent{\bf Step 5.}  Let us choose $\rho$, $\beta$ and $c$ such that the following holds
		\begin{align}\label{eq:first:order:condition}
		&C_1 = \beta- \rho A^T A \succ 0,\quad {\frac{c}{2} C_1 - \frac{2 C_1^T C_1}{\rho\tilde{\sigma}_{\min}(A^T A)}}\succ 0, \quad \frac{\beta}{2}- {cL}- \frac{2L^2}{\rho\tilde{\sigma}_{\min}(A^T A)}>0.
		\end{align}
		Note that the above inequalities are consistent, meaning there exists a tuple $(c, \beta, \rho)$ such that they will be satisfied simultaneously. One particular choice is
		\begin{align}
		\rho=\frac{16 L}{\tilde{\sigma}_{\min}(A^T A) }, \quad c =\frac{\beta}{2 L}
		\end{align}
		with $\beta$ chosen large enough such that
		\begin{align}
		\beta- \sigma_{\max}\left(\rho A^T A\right) -L > 0.
		\end{align}
		%These conditions can be satisfied by choosing $\theta$ satisfying \eqref{eq:theta}, $\rho$ being fixed, and $\beta$ to be large enough.
		
		%\noindent{\bf Step 6.} The boundedness of the $\lambda$ iteration can be see
		
		Then following similar argument as in Theorem 3.1 of \cite{hong16decomposing}, we can show that the first part of Claim \ref{claim:convergence} is true. {In particular, the boundedness of the primal and dual variable follows from part (2) of \cite{hong16decomposing} of Theorem 3.1, which utilizes Assumption [A1] and [A4].}  Further, by utilizing the standard argument in Theorem 2.9 of \cite{Bolte14}, we can claim the global convergence of the sequence $\{x^{r+1},\lambda^{r+1}\}$ under the $K{\L}$ assumption of $L(x,\lambda)$. We refer the readers to \cite{Li14splitting} for a similar argument.
		
		\subsection{Proof of Claim \ref{claim:sigma:T}}
		\noindent{\bf Proof.} First consider the matrix $Q^{-1}T$ defined in \eqref{eq:perturb}, which is given by
		\begin{align*}
		Q^{-1}T:= \begin{bmatrix}
		I_N & 0_{N\times M}\\ -\rho A & I_M
		\end{bmatrix}\begin{bmatrix}
		I -  \frac{1}{\beta}\left(H + \rho A^T A \right)& -\frac{1}{\beta} A^T\\ 0_{M\times N} & I_M
		\end{bmatrix}.
		\end{align*}
		Consider the characteristic polynomial of $Q^{-1}T$, given below
		\begin{align}\label{eq:determinant}
		\det (Q^{-1}T-\mu I) & =\det(T - \mu Q)\nonumber\\
		&=\det \begin{bmatrix}I_N -  \frac{1}{\beta}\left(H + \rho A^T A \right) -\mu I_N& -\frac{1}{\beta} A^T\\ \mu\rho A & I_M-\mu I_M
		\end{bmatrix}
		\end{align}
		where the first equality comes from the fact that $\det Q = 1$.
		First let us plug $\mu=1$ into the above equation.  We obtain
		\begin{align*}
		\det [T- Q]  &= \det \begin{bmatrix}I_N -  \frac{1}{\beta}\left(H + \rho A^T A \right) - I_N & -\frac{1}{\beta} A^T\\ \mu\rho A & 0
		\end{bmatrix} =0.
		\end{align*}
		Therefore we conclude that $\mu=1$ is an eigenvalue of $Q^{-1}T$.
		
		Second, let us test whether $0$ is an eigenvalue for $T$. To this end, plug $\mu=0$ into \eqref{eq:determinant}, we obtain
		\begin{align*}
		\det [T-0 \times Q]  &= \det \begin{bmatrix}I_N -  \frac{1}{\beta}\left(H + \rho A^T A \right) & -\frac{1}{\beta} A^T\\ 0_{M\times N}& I_M
		\end{bmatrix}.
		\end{align*}
		Therefore, as long as $\beta$ is {\it large enough} such that
		\begin{align}\label{eq:beta:1}
		\beta > \sigma_{\max} (H + \rho A^T A),
		\end{align}
		then the right hand side will not equal to zero. This suggests that  the matrix $Q^{-1}T$ is invertible if \eqref{eq:beta:1} holds true.
		
		Finally, let us investigate whether $Q^{-1}T$ has an eigenvalue which is strictly greater than $1$. For some $\delta>0$, let us take $\mu = 1+\delta$ and plug in to \eqref{eq:determinant}. We obtain
		\begin{align}
		\det [T-(1+\delta) Q]  &= \det \begin{bmatrix}I_N -  \frac{1}{\beta}\left(H + \rho A^T A \right) -(1+\delta) I_n & -\frac{1}{\beta} A^T\\ (1+\delta)\rho A & I_M-(1+\delta) I_M
		\end{bmatrix} \nonumber\\
		&=\det \left[-\delta (I_N -  \frac{1}{\beta}\left(H + \rho A^T A \right) -(1+\delta) I_N) + \frac{(1+\delta)\rho}{\beta}A^T A\right] \nonumber\\
		& = \det \left[\frac{\delta}{\beta}\left(H + \rho A^T A \right) + \delta^2 I_N + \frac{(1+\delta)\rho}{\beta}A^T A\right] \nonumber\\
		& = \det \left[\delta^2 I_N + \frac{\delta}{\beta} H + \frac{\rho (1 + 2\delta)}{\beta}A^T A\right] \nonumber.
		\end{align}
		We note that the matrix in the above determinant is symmetric therefore all its eigenvalues are real. To show that there exists a real $\delta>0$ such that the above determinant evaluates to zero, we follow the following two steps.
		
		\noindent{\bf Step 1.} It is clear that when $\delta>0$ is large enough, the $\delta^2$ term will dominate, and we have
		$$\delta^2 I_N + \frac{\rho}{\beta} (1+ 2 \delta)  A^T A + \frac{\delta}{\beta} H  \succ 0.$$
		Therefore the above determinant is positive.  Further, the matrix above will {\it only} have positive eigenvalues.
		
		\noindent{\bf Step 2.} Let us consider the case where $\delta$ is small and close to zero.
		Let $x^*$ be a strict saddle, then from \eqref{eq:strict} we have that there exists $y\in\mathbb{R}^n$ and $\sigma>0$ such that
		\begin{align}\label{eq:tilde:y}
		\|y\|=1, \; y\in \mbox{Null}(A), \; \mbox{and}\; y^T H y = -\sigma<0.
		\end{align}

		Note that the following holds true
		\begin{align*}
		&\delta^2 I_N+ \frac{\rho}{\beta}(1+ 2 \delta)  A^T A  + \frac{\delta}{\beta} H =\frac{1}{\beta} \left(\delta (H +\beta \delta I_N + 2 \rho A^T A) + \rho A^T A \right).
		\end{align*}
		It follows that for any $\rho>0$, if we choose $\delta = \sigma/2\beta$, then the following is true %matrix has {\it at least one} negative eigenvalue, which takes the value
		%$\frac{|\rho\sigma|\rho \sigma}{4\beta}$:
		\begin{align*}
		y^T\delta (H +\beta \delta I_N + 2 \rho A^T A) y =  y^T\delta (H +\beta \delta I_N ) y  = \delta(-\sigma  + \beta \delta) = -\delta \sigma/2<0.
		\end{align*}
		%This fact combined with the fact that $H+\gamma A^T A$ has at least one negative eigenvalue for any $\gamma>0$, implies that there exists $\tilde{y}\in\mbox{Null}(A)$ such that for any $\delta$ satisfying $\delta \le \frac{|\sigma(\rho)|}{2\beta}$, the following is true
		%\begin{align}\label{eq:tilde:y}
		%\delta \tilde{y}^T  (H +\beta \delta I) \tilde{y}<0.
		%\end{align}
		%To argue this, we define the projection matrix $P_{A} = A (A)^{+}$ where $A^{+}$ is the  Moore–Penrose pseudoinverse. Let $y\ne 0$, and $\|y\|=1$ be the vector such that
		%\begin{align}
		%y^T \delta (H +\beta \delta I + 2 \rho A^T A) y = \frac{|\sigma(\rho)|\sigma(\rho)}{4\beta}.
		%\end{align}
		%${y} =  (I-P_{A}) {y}' + P_{A} \hat{y}$ for some ${y}'$ and $\hat{y}$. Therefore, we have
		%\begin{align}
		%0>\frac{|\sigma(\rho)|\sigma(\rho)}{4\beta} & = y^T \delta (H +\beta \delta I + 2 \rho A^T A) y \nonumber\\
		%& = y^T \delta (H +\beta \delta I) y + \hat{y} (\delta (H +\beta \delta I)+ 2\rho A^T A) \hat{y}\nonumber\\
		%&\ge
		%\end{align}
		%
		%%\noindent{\bf Case I}. If for all $\tilde{y}$ such that $\tilde{y}^T A^T A \tilde{y} = 0$, we have $\tilde{y}^T H \tilde{y}= 0$. However, this case cannot happen because otherwise there exists $\bar{\gamma}$ large enough such that $H+ \gamma A^T A\succeq 0$, which contradicts to the strict saddle assumption \eqref{eq:strict}.
		%
		%\noindent{\bf Case II}. There exists $\tilde{y}$ such that $\tilde{y}^T H \tilde{y}\ne 0$, and $\tilde{y}^T A^T A \tilde{y} = 0$.

		Then we argue that for the same $\delta$, the following matrix also has at least one negative eigenvalue
		$$ \delta (H +\beta \delta I_N + 2 \rho A^T A) + \rho A^T A .$$
		
		%Suppose the above claim does not hold, then we have
		%\begin{align}\label{eq:contradiction}
		%z^T \left(\delta (H +\beta \delta I + 2 \rho A^T A) + \rho A^T A \right) z \ge 0, \; \forall~z\in\mathbb{R}^n.
		%\end{align}
		Use the ${y}$ given in \eqref{eq:tilde:y}, we have
		\begin{align*}
		{y}^T \left(\delta (H +\beta \delta I_N + 2 \rho A^T A) + \rho A^T A \right) {y} = \delta {y}^T  (H +\beta \delta I_N) {y} = -\delta\sigma/2<0.
		\end{align*}
		%which is a contradiction to \eqref{eq:contradiction}.
		
		Therefore we conclude that for any given $\beta>0$, there exists $\delta>0$ small enough, such that the following matrix has {\it at least one} negative eigenvalue
		\begin{align*}
		\delta^2 I_N + \frac{\rho}{\beta} (1+ 2 \delta)  A^T A + \frac{\delta}{\beta} H.
		\end{align*}
		By using the standard result on the continuity of eigenvalues, and by using the results in Step 1 and Step 2 above, we conclude that %for large enough $\delta>0$, the same matrix has {\it all} positive eigenvalues, we can apply  to conclude that
		there exists a positive $\delta^*>0$ such that
		\begin{align}
		\det (T-(1+\delta^*) Q)  & = \det\left((\delta^*)^2 I_N + \frac{\rho}{\beta} A^T A (1+ 2 \delta^*) + \frac{\delta^*}{\beta} H \right)=0. \label{eq:delta}
		\end{align}
		Let us refer to such an eigenvalue of $Q^{-1}T$ as $\sigma^*(Q^{-1}T)$, i.e., $\sigma^*(Q^{-1}T) = 1+\delta^*$.  \QED

		\subsection{Proof outline for Claim \ref{claim:convergence:admm}}
		In this section, we show the main steps leading to Claim \ref{claim:convergence:admm}.

		\noindent{\bf Step 1.} Let us define
		\begin{align*}
		&C_1 :=\beta I -\rho A^T A \succ 0, \quad  C_2: = \beta I -  \rho B^T B \succ 0.  \\
		&w^{r+1}:=(x^{r+1}-x^r)-(x^r-x^{r-1}) \quad v^{r+1}:=(y^{r+1}-y^r)-(y^r-y^{r-1}), \quad z^{r+1} = [w^{r+1}; v^{r+1}] . \\
		& W:=\begin{bmatrix}
		A^T A & A^T B \\ 0 & B^T B
		\end{bmatrix}, \quad  V:=[A, B], \quad C:=\beta I -\rho W.
		\end{align*}
		From the optimality condition of \eqref{eq:x:update:admm}- \eqref{eq:y:update:admm}, we have
		\begin{align*}
		\nabla f(x^r) + A^T \lambda^r + \rho A^T(Ax^{r+1} + B y^{r+1}-b) + \rho A^T B (y^r-y^{r+1}) + \rho A^T A (x^r-x^{r+1})+ \beta(x^{r+1}-x^r) =0,\\
		\nabla g(y^r) + B^T \lambda^r + \rho B^T(Ax^{r+1} + B y^{r+1}-b) + \rho B^T B (y^r-y^{r+1}) + \beta(y^{r+1}-y^r) =0.
		\end{align*}
		By some simple manipulation similarly as those leading to \eqref{eq:mu:difference:bound:1}, we can show
		\begin{align}
		\frac{1}{\rho}\|\lambda^{r+1}-\lambda^r\|^2\le \frac{2 L^2_g\|y^r-y^{r-1}\|^2 + 2L^2_f\|x^r-x^{r-1}\|}{\rho \tilde{\sigma}_{\min}(V^T V)} + \frac{2}{\rho \tilde{\sigma}_{\min}(V^T V)}\|z^{r+1}\|_{C^T C}^2
		\end{align}
		
		\noindent{\bf Step 2.} From the optimality conditions of \eqref{eq:x:update:admm} and \eqref{eq:y:update:admm}, we have
		\begin{align*}
		&\langle \nabla g(y^r) + B^T \lambda^{r+1} +(\beta I -\rho B^T B) (y^{r+1}-y^r), y-y^r\rangle =0, \; \forall~y \\
		&\langle  \nabla f(x^r) + A^T \lambda^{r+1} +(\beta I -\rho A^T A) (x^{r+1}-x^r) + \rho A^T B (y^r-y^{r+1}), x-x^r \rangle =0,\; \forall~x.
		\end{align*}
		Then subtracting the previous iteration of the same condition, and add them together, we obtain
		\begin{align}
		&\langle \lambda^{r+1}-\lambda^r, B(y^{r+1}-y^r) + A(x^{r+1}-x^r)\rangle \nonumber\\
		&\le \langle \nabla g(y^r)-\nabla g(y^{r-1}), y^r-y^{r+1}\rangle + \langle \nabla f(x^r)-\nabla f(x^{r-1}), x^r-x^{r+1}\rangle \nonumber\\
		&\quad\quad- \langle (\beta I -\rho B^T B) v^{r+1}, y^{r+1}-y^r\rangle -\langle (\beta I -\rho A^T A) w^{r+1}, x^{r+1}-x^r\rangle +\langle \rho A^T B v^{r+1}, x^{r+1}-x^r\rangle.
		\end{align}
		%Further, $\beta$ is chosen such that there exists $\xi>0$
		%\begin{align}
		%\beta I -\rho A^T A := C_1 \succeq   \xi I \quad \beta I -2 {\rho} B^T B :=C_2 \succeq  \xi I. \label{eq:AB:bound}
		%\end{align}
		Collecting terms, and after some simple manipulation, we obtain
		\begin{align}\label{eq:dual:descent}
		&\frac{1}{2\rho}\|\lambda^{r+1}-\lambda^r\|^2 + \frac{1}{2}\|x^{r+1}-x^r\|_{C_3}^2+ \frac{1}{2}\|y^{r+1}-y^r\|^2_{C_4}\nonumber\\
		&\le \frac{1}{2\rho} \|\lambda^{r}-\lambda^{r-1}\|^2+ \frac{1}{2}\|x^{r}-x^{r-1}\|_{C_3}^2+ \frac{1}{2}\|y^{r}-y^{r-1}\|_{C_4}^2 \nonumber\\
		& \quad +\frac{1}{2}\|x^{r+1}-x^r\|_{C_5}^2+\frac{1}{2} \|y^{r+1}-y^r\|_{C_6}^2 -\frac{1}{2}\|w^{r+1}\|_{C_1}^2 - \frac{1}{2}\|v^{r+1}\|_{C_2-\rho B^T B}^2.
		\end{align}
		where
		\begin{align}
		C_3 &= (\beta I -\rho A^T A)  + {L_f}I, \quad C_4 = (\beta I -\rho B^T B) +   {L_g}I, \\
		C_5 &= 2L_f + \rho A^T A ,\quad C_6 = 2L_g I.
		\end{align}
		
		\noindent{\bf Step 3.} By a standard descent estimate for the gradient-type algorithm, we can show that the augmented Lagrangian decreases in the following manner
		\begin{align}\label{eq:primal:descent}
		&L(x^{r+1}, y^{r+1}, \lambda^{r+1}) - L(x^{r}, y^{r}, \lambda^{r}) \nonumber\\
		&\le -\frac{\beta}{2}\|x^{r+1}-x^r\|^2 -\frac{\beta}{2}\|y^{r+1}-y^r\|^2 + \frac{1}{\rho}\|\lambda^{r+1}-\lambda^r\|^2\nonumber\\
		&\le -\frac{\beta}{2}\|x^{r+1}-x^r\|^2 -\frac{\beta}{2}\|y^{r+1}-y^r\|^2 + \frac{2 L^2_g\|y^r-y^{r-1}\|^2+ 2L^2_f\|x^r-x^{r-1}\|^2}{\rho \tilde{\sigma}_{\min}(V^T V)} + \frac{2}{\rho \tilde{\sigma}_{\min}(V^T V)}\|z^{r+1}\|_{C^T C}^2
		\end{align}
		whenever the following holds
		\begin{align}
		\beta- L_f -\sigma_{\max}(\rho A^T A)>0, \quad \beta- L_g -\sigma_{\max}(\rho B^T B)>0.
		\end{align}
		Adding the above inequality with \eqref{eq:dual:descent} multiplied by a constant $2 c>0$, we obtain
		\begin{align}
		&L(x^{r+1}, y^{r+1}, \lambda^{r+1}) + \frac{c}{\rho}\|\lambda^{r+1}-\lambda^r\|^2 + c\|x^{r+1}-x^r\|_{C_3}^2+ \left(\frac{2 L^2_f}{\rho \tilde{\sigma}_{\min}(V^T V)} \right)\|x^{r+1}-x^r\|^2 \nonumber\\
		&\quad +\left(\frac{2 L^2_g}{\rho \tilde{\sigma}_{\min}(V^T V)} \right)\|y^{r+1}-y^r\|^2 + c\|y^{r+1}-y^r\|^2_{C_4}\nonumber\\
		&\le   L(x^{r}, y^{r}, \lambda^{r}) + \frac{c}{\rho}\|\lambda^{r}-\lambda^{r-1}\|^2+ c\|x^{r}-x^{r-1}\|_{C_3}^2 + \left(\frac{2 L^2_f}{\rho \tilde{\sigma}_{\min}(V^T V)} \right)\|x^{r+1}-x^r\|^2 \nonumber\\
		&\quad + \left(\frac{2 L^2_g}{\rho \tilde{\sigma}_{\min}(V^T V)} \right)\|y^{r}-y^{r-1}\|^2 + c\|y^{r+1}-y^r\|^2_{C_4}\nonumber\\
		&\quad -\left(\frac{\beta}{2}-c(2L_f + \sigma_{\max}(A^T A) \rho)- \frac{2 L^2_f}{\rho \tilde{\sigma}_{\min}(V^T V)}\right)\|x^{r+1}-x^r\|^2-\left(\frac{\beta}{2}-2cL_g-\frac{2 L^2_g}{\rho \tilde{\sigma}_{\min}(V^T V)} \right)\|y^{r+1}-y^r\|^2\nonumber\\
		&\quad -(z^{r+1})^T\left(\frac{c}{2} \begin{bmatrix}
		C_1 & 0 \\ 0 & C_2-\rho B^T B
		\end{bmatrix} -\frac{2 C^T C}{\rho \tilde{\sigma}_{\min}(V^T V)}\right) z^{r+1}.
		\end{align}
		Therefore, to make the entire potential function decrease, we will need the following conditions
		\begin{align}\label{eq:condition:admm}
		\begin{split}
		& \beta  - \sigma_{\max}(\rho B^T B)- L_g > 0, \quad \beta - \sigma_{\max}(\rho A^T A)- L_f >0 \\
		&\frac{c}{2} \begin{bmatrix}
		C_1 & 0 \\ 0 & C_2-\rho B^T B
		\end{bmatrix} -\frac{2 C^T C}{\rho \tilde{\sigma}_{\min}(V^T V)}\succ 0, \\
		& \frac{\beta}{2}-2 c L_f -c \sigma_{\max}(A^T A)\rho - \frac{2 L^2_f}{\rho \tilde{\sigma}_{\min}(V^T V)}>0, \\
		& \frac{\beta}{2}-2 c L_g -\frac{2 L^2_g}{\rho \tilde{\sigma}_{\min}(V^T V)} >0.
		\end{split}
		\end{align}
		Similarly as argued in \eqref{eq:first:order:condition}, these inequalities are consistent, meaning there exists a choice of $\beta, c, \rho$ such that they will all be satisfied.
		
		The rest of the steps are similar to Step 3 - Step 5 in the outline of the proof of Claim \ref{claim:convergence}. We do not repeat them here.
		
		\subsection{Proof of Claim \ref{claim:real}}
		
		\noindent{\bf Proof.}
		By using standard determinant for block matrices, we obtain
		\begin{align*}
		&\det[T-\mu P] \nonumber\\
		&= \det\begin{bmatrix}
		(1-\mu) I_N -\frac{1}{\beta} H -\frac{\rho}{\beta}A^T A & -\frac{\rho}{\beta}A^T B & -\frac{1}{\beta}A^T\\
		-\mu\frac{\rho}{\beta}A^T A & (1-\mu) I_N -\frac{1}{\beta} G +\frac{\rho}{\beta}B^T B & -\frac{1}{\beta} B^T\\
		\rho\mu A & \rho \mu B & (1-\mu) I_M
		\end{bmatrix} \nonumber\\
		&= (1-\mu)\det\left(\begin{bmatrix}
		(1-\mu) I -\frac{1}{\beta} H -\frac{\rho}{\beta}A^T A & -\frac{\rho}{\beta}A^T B \\
		-\mu\frac{\rho}{\beta}B^T A & (1-\mu) I_N -\frac{1}{\beta} G +\frac{\rho}{\beta}B^T B \\
		\end{bmatrix} - \frac{1}{1-\mu} \begin{bmatrix}
		-\frac{\rho \mu}{\beta}A^T A & -\frac{\rho \mu}{\beta}A^T B \\
		-\frac{\rho \mu}{\beta}B^T A & -\frac{\rho\mu}{\beta}B^T B
		\end{bmatrix}\right)\\
		& := \det[U(\mu)]
		\end{align*}
		where we have defined the matrix $U(\mu) = [U_{11}(\mu)\; U_{12}(\mu); U_{12}(\mu)\; U_{22}(\mu)]\in\mathbb{R}^{2N\times 2N}$, with
		%\begin{align}
		%U_{11}(\mu)\nonumber\\
		%U_{12}(\mu)\nonumber\\
		%U_{21}(\mu)\nonumber\\
		%U_{22}(\mu)\nonumber\\
		%\end{align}
		%
		
		%\begin{align*}
		%\det\left[T-\mu M \right] & \stackrel{\rm(i)}= \det\left[-\mu(T_{11}-M_{11}\mu) - T_{12}T_{21} \right]\nonumber\\
		%&=\det\Bigg[-\mu\begin{bmatrix}
		%2I - \frac{2\rho}{\beta}A^T A-\frac{1}{\beta} H -\mu I & -\frac{2\rho}{\beta}A^T B \\
		%-\frac{\rho}{\beta}B^T A \mu & 2 I -\frac{1}{\beta}G -\frac{2\rho}{\beta}B^T B -\mu I
		%\end{bmatrix} \nonumber\\
		%&\quad - \begin{bmatrix}
		%-I + \frac{\rho}{\beta}A^T A +\frac{1}{\beta} H & \frac{\rho}{\beta}A^T B \\
		% 0 & -I+\frac{1}{\beta} G+ \frac{\rho}{\beta}B^T B
		%\end{bmatrix}\Bigg]\nonumber\\
		%& := \det[Q(\mu)]
		%\end{align*}
		%where $(i)$ comes from \cite[Theorem 3]{Silvester00} on  the determinant of $2\times 2$ block matrices with two of the blocks commute; we have defined the matrix $Q(\mu) = [Q_{11}(\mu)\; Q_{12}(\mu); Q_{12}(\mu)\; Q_{22}(\mu)]\in\mathbb{R}^{2N\times 2N}$, with
		\begin{subequations}
			\begin{align}
			U_{11}(\mu) &= -\mu\left(2I - \frac{2\rho}{\beta}A^T A-\frac{1}{\beta} H -\mu I\right) + I - \frac{\rho}{\beta}A^T A -\frac{1}{\beta} H\\
			U_{12}(\mu) & = \mu \frac{2\rho}{\beta}A^T B -\frac{\rho}{\beta}A^T B = (2\mu-1)\frac{\rho}{\beta}A^T B \\
			U_{21}(\mu) & = \mu^2\frac{\rho}{\beta}B^T A\\
			U_{22}(\mu) & = -\mu\left(2 I -\frac{1}{\beta}G -\frac{2\rho}{\beta}B^T B -\mu I\right) + I-\frac{1}{\beta} G- \frac{\rho}{\beta}B^T B.
			\end{align}
		\end{subequations}
		We first verify the case with $\mu=1$. In this case, it is easy to verify that
		\begin{align*}
		U (1)= \frac{\rho}{\beta}\begin{bmatrix}
		A^T A & A^T B \\ B^T A & B^T B
		\end{bmatrix}
		\end{align*}
		therefore $\det[U(1)]=0$, implying that $\mu=1$ is an eigenvalue for the matrix $P^{-1} T$.
		
		Also let $\mu=0$, we have
		\begin{align*}
		U(0) = \begin{bmatrix}
		I - \frac{\rho}{\beta}A^T A -\frac{1}{\beta} H & -\frac{\rho}{\beta}A^T B\\
		0 & I-\frac{1}{\beta} G- \frac{\rho}{\beta}B^T B
		\end{bmatrix}.
		\end{align*}
		Clearly, when $\beta$ satisfies the following inequalities, the matrix is invertible
		\begin{align*}
		\beta> \rho\sigma_{\max}(A^T A) + L_f, \quad \beta> \rho\sigma_{\max}(B^T B) + L_g,
		\end{align*}
		where $L_f$ and $L_g$ represent the Lipschitz constant for the objective function $\nabla f$ and $\nabla g$.

		We note that $U(1+\delta)$ can be written in the following form
		\begin{align}\label{eq:key}
		U(1+\delta)&=\begin{bmatrix}
		U_{11}(1+\delta) & (2\delta+1) \frac{\rho}{\beta} A^T B\\ (2\delta+1+\delta^2) \frac{\rho}{\beta} B^T A & U_{22}(1+\delta)
		\end{bmatrix} \nonumber\\
		&=\begin{bmatrix}
		U_{11}(1+\delta) & (2\delta+1) \frac{\rho}{\beta}A^T B\\  \frac{(\delta+1)^2 }{2\delta+1} (2\delta+ 1) \frac{\rho}{\beta}B^T A & U_{22}(1+\delta)
		\end{bmatrix} \nonumber\\
		&=\begin{bmatrix}
		I & 0 \\ 0 & \frac{\delta+1}{\sqrt{2\delta+1}}
		\end{bmatrix}
		\begin{bmatrix}
		U_{11}(1+\delta) & {(\delta+1)}{\sqrt{2\delta+1}} \frac{\rho}{\beta}A^T B\\ {(\delta+1)}{\sqrt{2\delta+1}} \frac{\rho}{\beta}B^T A & U_{22}(1+\delta)
		\end{bmatrix} \begin{bmatrix}
		I  & 0 \\ 0 & \frac{\sqrt{2\delta+1}}{\delta+1}
		\end{bmatrix}
		\nonumber.
		\end{align}
		By noting the fact that $\delta+1>0$, and $U_{11}(1+\delta)$ and $U_{22}(1+\delta)$ are both symmetric matrices, we conclude that $U(1+\delta)$ has real eigenvalues. \QED
		
		\section{Numerical Results}
		
\begin{figure*}[ht]
			\vspace{-0.3cm}
			\begin{minipage}[t]{0.48\linewidth}
				\centering
				{\includegraphics[width=0.75\linewidth]{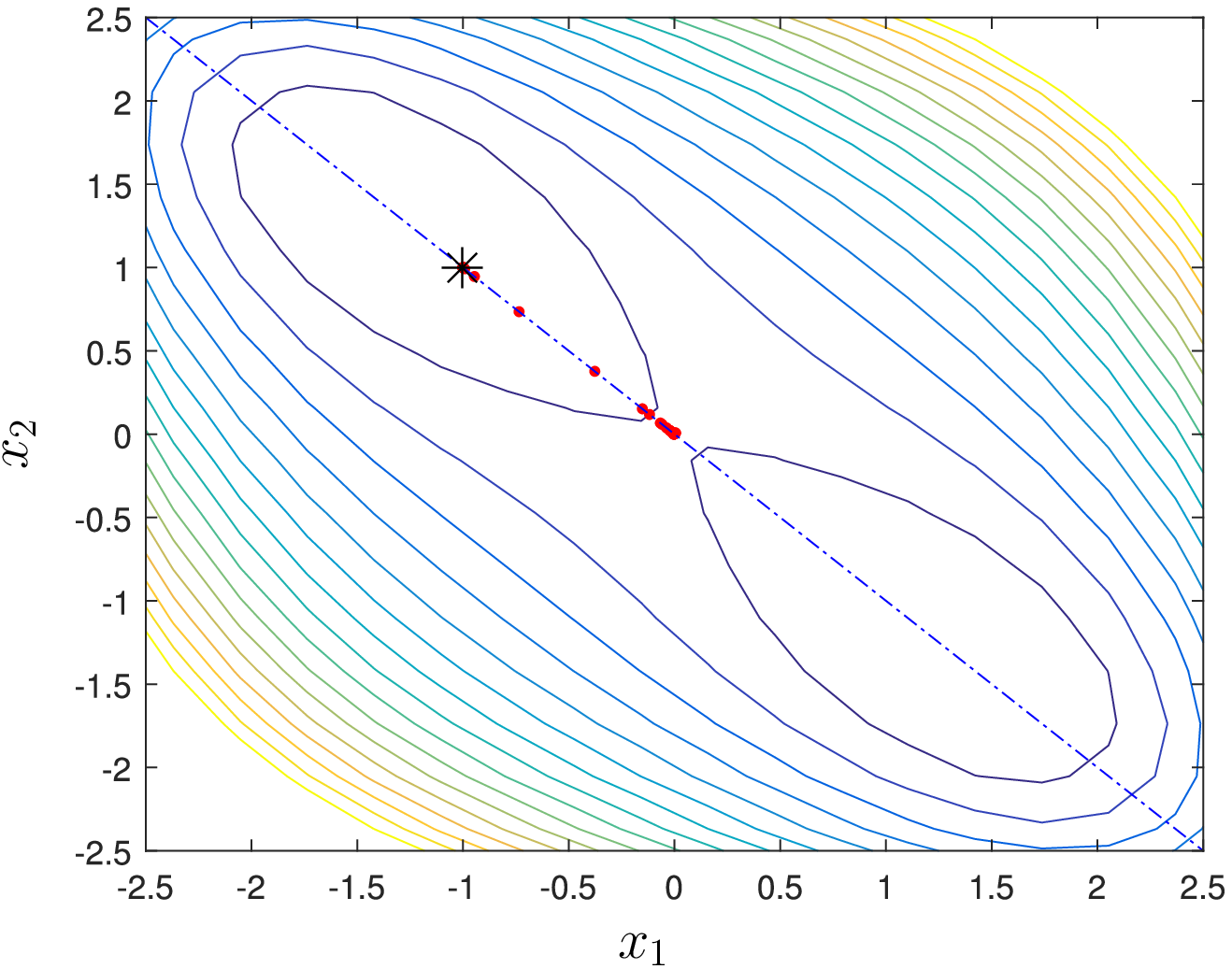}}%\label{figcou}
			\end{minipage}
			\hfill
			\begin{minipage}[t]{0.48\linewidth}
				{\includegraphics[width=0.75\linewidth]{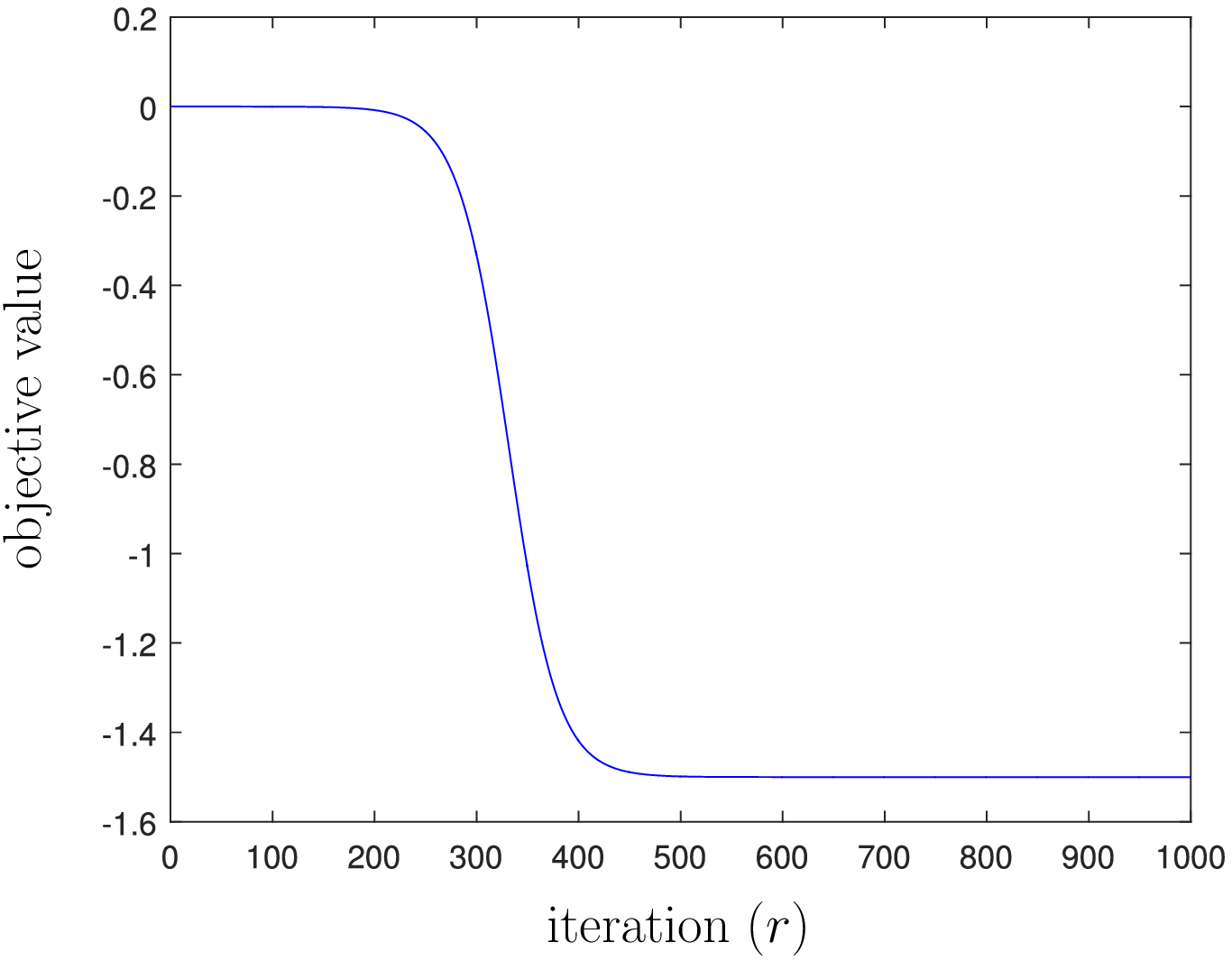} }%\label{figobj2}
			\end{minipage}
			\vspace{-0.2cm}
			\caption{Contour of the objective function and trajectory of the iterates, where $N=2$, $M=1$, $\rho=10$, $\beta=200$.}
			\label{fig:con}
		\end{figure*}

		Consider a nonconvex objective function
		\begin{equation}
		f(x)=x^{T}Qx+\frac{1}{4}\|x\|^4_4\label{eq.siobj2}
		\end{equation}
		where $Q\in\mathbb{R}^{N\times N}$ is indefinite.
		First, we have the following properties of function $f(x)$ such that $f(x)$ satisfies the assumptions of the analysis.
		\begin{lemma}\label{le.simup}
			For any $\tau\ge\lambda_{\max}(Q)$ and $x\in\{x|\|x\|^2\le\tau\}$, $f(x)$ defined in \eqref{eq.siobj2} is $5\tau$-smooth and $6\sqrt{\tau}$-Hessian Lipschitz.
		\end{lemma}
		\nonumber{\bf 2-D Case.} We first test GPDA on a two dimensional case, where $A=[1 -1]$ and $b=1$. Constraint $Ax=b$ forms a line in this case which is shown in Fig. \ref{fig:con} with blue color. The GPDA algorithm is randomly initialized at the origin. It can be seen from Fig.  \ref{fig:con}  that there are two local optimal points and one strict saddle point at the origin, and GPDA can escape from the saddle point efficiently. Since there is a constraint, the iterates converge to a point on the line which is the nearest point to the local optimal point of the objective function. 

		\paragraph{Random matrix $Q$.} We also randomly generate matrix $Q$ with the following steps: 
		
		\noindent 1) randomly generate a diagonal matrix $D$ whose entries follow \emph{i.i.d.} Gaussian distribution with zero mean and variance one; 
		
		\noindent 2) generate an orthogonal matrix $U\in\mathbb{R}^{d\times d}$; 
		
		\noindent  3) obtain matrix $Q=UDU^{T}$. The entries of matrix $A\in\mathbb{R}^{M\times N}$ where $M=5$ and $b$ are also randomly generated, which follow \emph{i.i.d.} Gaussian distribution with zero mean and variance one. 
		
\begin{figure*}[ht]
			\vspace{-0.3cm}
			\begin{minipage}[t]{0.48\linewidth}
				\centering
				{\includegraphics[width=0.75\linewidth]{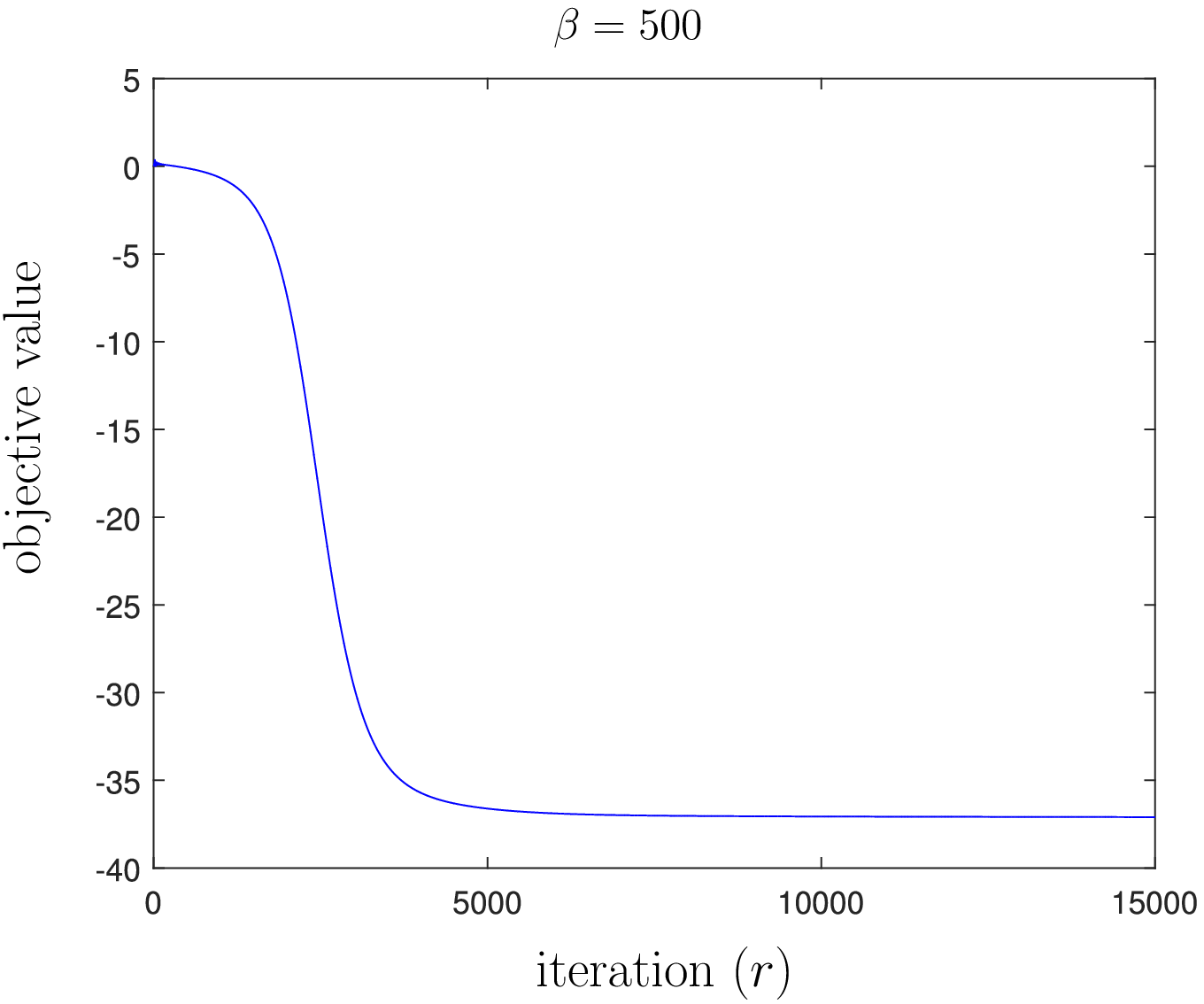} }	\label{fig:cc}
			\end{minipage}
			\hfill
			\begin{minipage}[t]{0.48\linewidth}
				{\includegraphics[width=	0.75\linewidth]{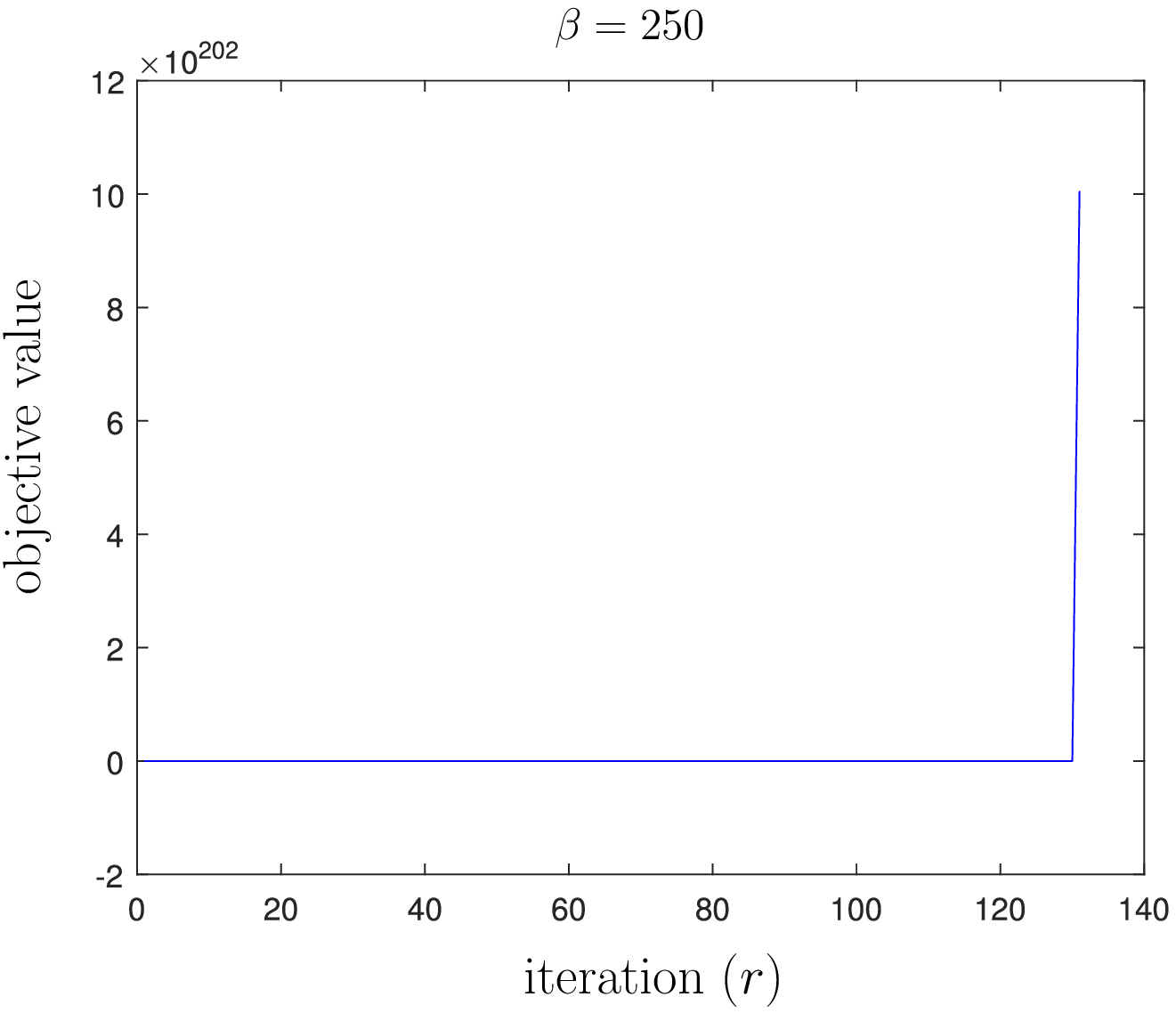}  }	\label{fig:nc}
			\end{minipage}
			\vspace{-0.2cm}
			\caption{Objective value of function \eqref{eq.siobj2}, where $N=20$, $M=5$, $\rho=10$}
			\label{fig:con:2}
		\end{figure*}

		We initialize GPDA around the strict saddle point which is at the origin randomly. It can be observed from Fig. \ref{fig:con:2} that GPDA can converge to a point where the corresponding objective value is much lower than the one at the origin, implying that GPDA can escape from the saddle point efficiently. Also, we can see that when $\beta$ is small, GPDA will diverge, indicating that $\beta$ should be large enough, which is consistent with the theoretical analysis.
		
		\subsection{Numerical Example}%\leref{le.simup}
		\begin{proof}
			Consider function
			\begin{equation}
			f(x)=\bx^{T}Q\bx+\frac{1}{4}\|x\|^4_4\label{eq.siobj2:2}
			\end{equation}
			where $x\in\mathcal{S}$, $\mathcal{S}=\{x|\|x\|^2\le\tau\}$ and $\tau\ge\lambda_{\max}(Q)$.
			
			\paragraph{To prove L-smooth Lipschitz continuity}:
			\begin{align}
			\notag
			\|\nabla f(x)-\nabla f(y)\|=&\left\|2(Qx-Qy)+\left[\begin{array}{c}x^3_1-y^3_1\\\vdots\\x^3_d-y^3_d\end{array}\right]\right\|,\quad\forall x,y\in\mathcal{S}
			\\\notag
			\le&2\lambda_{\max}(Q)\|x-y\|+\left\|\left[\begin{array}{c}(x_1-y_1)(x^2_1+x_1y_1+y^2_1)\\\vdots\\(x_d-y_d)(x^2_d+x_dy_d+y^2_d)\end{array}\right]\right\|
			\\\notag
			\mathop{\le}\limits^{(a)}&2\tau\|x-y\|+3\tau\|x-y\|\le5\tau\|x-y\|
			\end{align}
			where $x_i$ denotes the $i$th entry of vector $\bx$, and $(a)$ is true because
			\begin{equation}
			x^2_i\le\tau,\quad y^2_i\le\tau,\quad x_iy_i\le (x^2_i+y^2_i)/2\le\tau,\forall i.\label{eq.scain}
			\end{equation}
			
			\paragraph{To prove Hessian Lipschitz continuity}:
			\begin{align}\notag
			\|\nabla^2 f(x)-\nabla^2 f(y)\|=&3\left\|\begin{array}{ccc}x^2_1-y^2_1 & \cdots & 0 \\ \vdots & \ddots & \vdots \\ 0 & \cdots & x^2_d-y^2_d\end{array}\right\|
			\\\notag
			\le&6\sqrt{\tau}\left\|\begin{array}{ccc}x_1-y_1 & \cdots & 0 \\ \vdots & \ddots & \vdots \\ 0 & \cdots & x_d-y_d\end{array}\right\|=6\sqrt{\tau}\|x-y\|
			\end{align}
			where $(a)$ is true because $x_i+y_i\le\sqrt{(x_i+y_i)^2}=\sqrt{x^2_1+2x_iy_i+y^2_i}\mathop{\le}\limits^{\eqref{eq.scain}}2\sqrt{\tau},\forall i$.
			
		\end{proof}

{\small
\bibliographystyle{IEEEtran}
\bibliography{ref,biblio,ref_pprox_pda,refs,ref_zero1,ref_zero2}}

\end{document}